# PINNING AND WETTING TRANSITION FOR (1+1)-DIMENSIONAL FIELDS WITH LAPLACIAN INTERACTION

By Francesco Caravenna[1] and Jean-Dominique Deuschel

*Università degli Studi di Padova and Technische Universität Berlin*

We consider a random field $\varphi : \{1, \ldots, N\} \to \mathbb{R}$ as a model for a linear chain attracted to the defect line $\varphi = 0$, that is, the $x$-axis. The free law of the field is specified by the density $\exp(-\sum_i V(\Delta\varphi_i))$ with respect to the Lebesgue measure on $\mathbb{R}^N$, where $\Delta$ is the discrete Laplacian and we allow for a very large class of potentials $V(\cdot)$. The interaction with the defect line is introduced by giving the field a reward $\varepsilon \geq 0$ each time it touches the $x$-axis. We call this model the *pinning model*. We consider a second model, the *wetting model*, in which, in addition to the pinning reward, the field is also constrained to stay nonnegative.

We show that both models undergo a phase transition as the intensity $\varepsilon$ of the pinning reward varies: both in the pinning ($a = \mathrm{p}$) and in the wetting ($a = \mathrm{w}$) case, there exists a critical value $\varepsilon_c^a$ such that when $\varepsilon > \varepsilon_c^a$ the field touches the defect line a positive fraction of times (localization), while this does not happen for $\varepsilon < \varepsilon_c^a$ (delocalization). The two critical values are nontrivial and distinct: $0 < \varepsilon_c^{\mathrm{p}} < \varepsilon_c^{\mathrm{w}} < \infty$, and they are the only nonanalyticity points of the respective free energies. For the pinning model the transition is of second order, hence the field at $\varepsilon = \varepsilon_c^{\mathrm{p}}$ is delocalized. On the other hand, the transition in the wetting model is of *first order* and for $\varepsilon = \varepsilon_c^{\mathrm{w}}$ the field is localized. The core of our approach is a Markov renewal theory description of the field.

## 1. Introduction and main results.

1.1. *Definition of the models.* We are going to define two distinct but related models for a (1+1)-dimensional random field. These models depend

Received March 2007; revised October 2007.

[1]Supported in part by the German Research Foundation–Research Group 718 during his stay at TU Berlin in January 2007.

*AMS 2000 subject classifications.* 60K35, 60F05, 82B41.

*Key words and phrases.* Pinning model, wetting model, phase transition, entropic repulsion, Markov renewal theory, local limit theorem, Perron–Frobenius theorem.

on a measurable function $V(\cdot):\mathbb{R}\to\mathbb{R}\cup\{+\infty\}$, the *potential*. We require that $x\mapsto\exp(-V(x))$ is bounded and continuous and that $\int_{\mathbb{R}}\exp(-V(x))\,dx<\infty$. Since a global shift on $V(\cdot)$ is irrelevant for our purposes, we will actually impose the stronger condition

$$\int_{\mathbb{R}} e^{-V(x)}\,dx = 1. \tag{1.1}$$

The last assumptions we make on $V(\cdot)$ are that $V(0)<\infty$, that is, $\exp(-V(0))>0$, and that

$$\int_{\mathbb{R}} x^2 e^{-V(x)}\,dx =: \sigma^2 < \infty \quad \text{and} \quad \int_{\mathbb{R}} x\, e^{-V(x)}\,dx = 0. \tag{1.2}$$

A typical example is of course the Gaussian case $V(x)=x^2/(2\sigma^2)+\log\sqrt{2\pi\sigma^2}$, but we stress that we do not make any convexity assumption on $V(\cdot)$. Next we introduce the Hamiltonian $\mathcal{H}_{[a,b]}(\varphi)$, defined for $a,b\in\mathbb{Z}$, with $b-a\geq 2$, and for $\varphi:\{a,\ldots,b\}\to\mathbb{R}$ by

$$\mathcal{H}_{[a,b]}(\varphi) := \sum_{n=a+1}^{b-1} V(\Delta\varphi_n), \tag{1.3}$$

where $\Delta$ denotes the discrete Laplacian:

$$\Delta\varphi_n := (\varphi_{n+1}-\varphi_n)-(\varphi_n-\varphi_{n-1}) = \varphi_{n+1}+\varphi_{n-1}-2\varphi_n. \tag{1.4}$$

We are ready to introduce our first model, the *pinning model* (p-*model* for short) $\mathbb{P}^{\mathrm{p}}_{\varepsilon,N}$, that is, the probability measure on $\mathbb{R}^{N-1}$ defined by

$$\mathbb{P}^{\mathrm{p}}_{\varepsilon,N}(d\varphi_1\cdots d\varphi_{N-1}) := \frac{\exp(-\mathcal{H}_{[-1,N+1]}(\varphi))}{\mathcal{Z}^{\mathrm{p}}_{\varepsilon,N}} \prod_{i=1}^{N-1}(\varepsilon\delta_0(d\varphi_i)+d\varphi_i) \tag{1.5}$$

where $N\in\mathbb{N}$, $\varepsilon\geq 0$, $d\varphi_i$ is the Lebesgue measure on $\mathbb{R}$, $\delta_0(\cdot)$ is the Dirac mass at zero and $\mathcal{Z}^{\mathrm{p}}_{\varepsilon,N}$ is the normalization constant, usually called partition function. To complete the definition, in order to make sense of $\mathcal{H}_{[-1,N+1]}(\varphi)$, we have to specify:

$$\text{the boundary conditions} \quad \varphi_{-1}=\varphi_0=\varphi_N=\varphi_{N+1}:=0. \tag{1.6}$$

We fix zero boundary conditions for simplicity, but our approach works for arbitrary choices (as long as they are bounded in $N$).

The second model we consider, the *wetting model* (w-*model* for short) $\mathbb{P}^{\mathrm{w}}_{\varepsilon,N}$, is a variant of the pinning model defined by

$$\begin{aligned}\mathbb{P}^{\mathrm{w}}_{\varepsilon,N}(d\varphi_1\cdots d\varphi_{N-1}) &:= \mathbb{P}^{\mathrm{p}}_{\varepsilon,N}(d\varphi_1\cdots d\varphi_{N-1}\,|\,\varphi_1\geq 0,\ldots,\varphi_{N-1}\geq 0)\\ &= \frac{\exp(-\mathcal{H}_{[-1,N+1]}(\varphi))}{\mathcal{Z}^{\mathrm{w}}_{\varepsilon,N}} \prod_{i=1}^{N-1}(\varepsilon\delta_0(d\varphi_i)+d\varphi_i\mathbf{1}_{(\varphi_i\geq 0)}),\end{aligned} \tag{1.7}$$

that is, we replace the measure $d\varphi_i$ by $d\varphi_i \mathbf{1}_{(\varphi_i \geq 0)}$ and $\mathcal{Z}^{\mathrm{p}}_{\varepsilon,N}$ by a new normalization $\mathcal{Z}^{\mathrm{w}}_{\varepsilon,N}$.

Both $\mathbb{P}^{\mathrm{p}}_{\varepsilon,N}$ and $\mathbb{P}^{\mathrm{w}}_{\varepsilon,N}$ are (1+1)-dimensional models for a linear chain of length $N$ which is attracted to a defect line, the $x$-axis, and the parameter $\varepsilon \geq 0$ tunes the strength of the attraction. By "(1+1)-dimensional" we mean that the configurations of the linear chain are described by the trajectories $\{(i, \varphi_i)\}_{0 \leq i \leq N}$ of the field, so that we are dealing with directed models (see Figure 1 for a graphical representation). We point out that linear chain models with Laplacian interaction appear naturally in the physical literature in the context of *semiflexible polymers*; cf. [6, 17, 21] (however, the scaling they consider is different from the one we look at in this paper). One note about the terminology: while "pinning" refers of course to the attraction terms $\varepsilon\delta_0(d\varphi_i)$, the use of the term "wetting" is somewhat customary in the presence of a positivity constraint and refers to the interpretation of the field as an effective model for the interface of separation between a liquid above a wall and a gas; see [10] for more details.

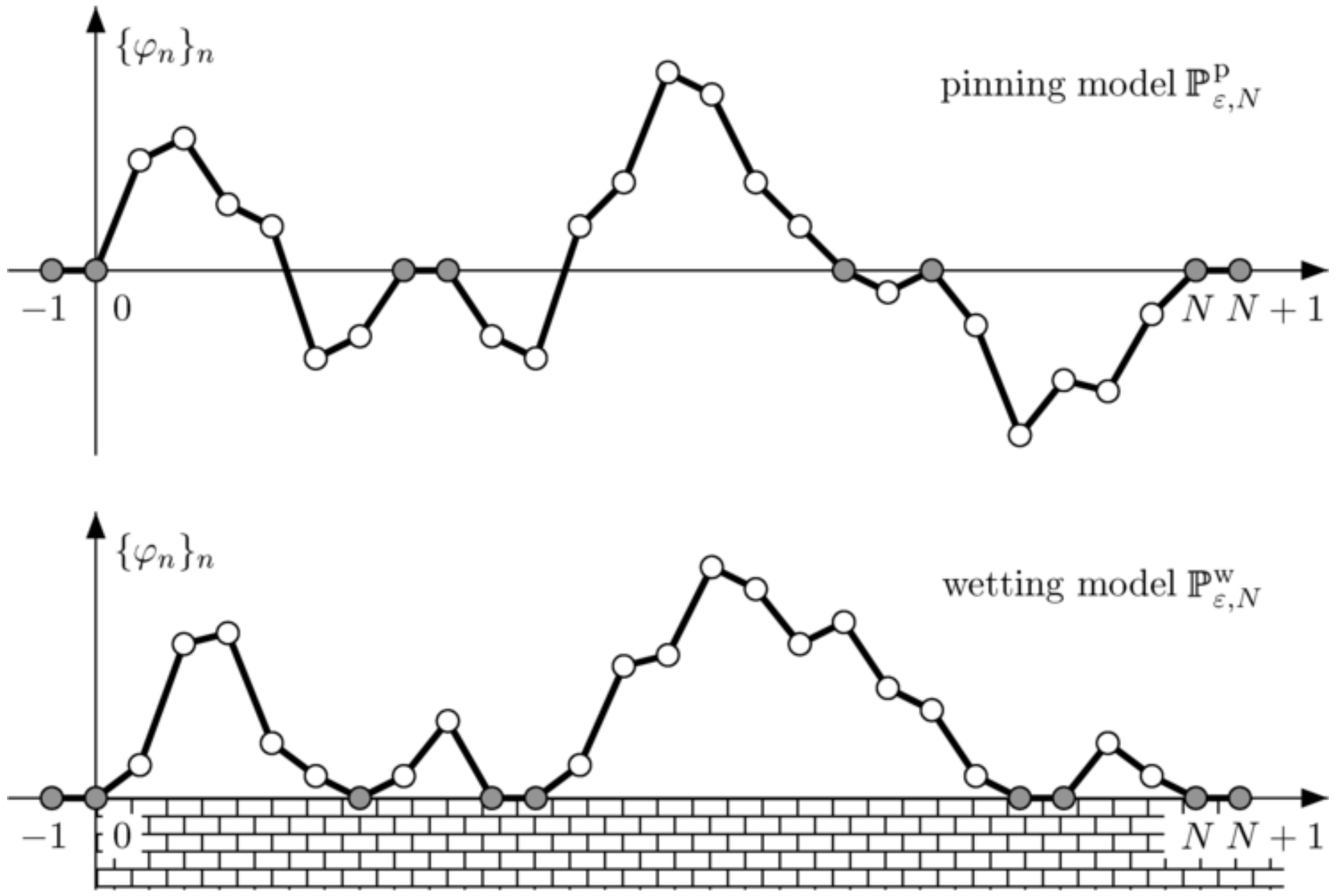

FIG. 1. *A graphical representation of the pinning model $\mathbb{P}^{\mathrm{p}}_{\varepsilon,N}$ (top) and of the wetting model $\mathbb{P}^{\mathrm{w}}_{\varepsilon,N}$ (bottom), for $N = 25$ and $\varepsilon > 0$. The trajectories $\{(n, \varphi_n)\}_{0 \leq n \leq N}$ of the field describe the configurations of a linear chain attracted to a defect line, the $x$-axis. The gray circles represent the pinned sites, that is, the points in which the chain touches the defect line, which are energetically favored. Note that in the pinning case the chain can cross the defect line without touching it, while this does not happen in the wetting case due to the presence of a wall, that is, of a constraint for the chain to stay nonnegative; the repulsion effect of entropic nature that arises is responsible for the different critical behavior of the models.*

The purpose of this paper is to investigate the behavior of $\mathbb{P}^{\mathrm{p}}_{\varepsilon,N}$ and $\mathbb{P}^{\mathrm{w}}_{\varepsilon,N}$ in the large $N$ limit; in particular we wish to understand whether and when the reward $\varepsilon \geq 0$ is strong enough to pin the chain at the defect line, a phenomenon that we will call *localization*. We point out that these kinds of questions have been answered in depth in the case of gradient interaction, that is, when the Laplacian $\Delta$ appearing in (1.3) is replaced by the discrete gradient $\nabla\varphi_n := \varphi_n - \varphi_{n-1}$ (cf. [2, 9, 10, 11, 13, 18]); we will refer to this as the *gradient case*. As we are going to see, the behavior in the Laplacian case turns out to be sensibly different.

1.2. *The free energy and the main results.* A convenient way to define localization for our models is by looking at the Laplace asymptotic behavior of the partition function $\mathcal{Z}^a_{\varepsilon,N}$ as $N \to \infty$. More precisely, for $a \in \{\mathrm{p}, \mathrm{w}\}$ we define the *free energy* $\mathrm{F}^a(\varepsilon)$ by

$$\mathrm{F}^a(\varepsilon) := \lim_{N\to\infty} \mathrm{F}^a_N(\varepsilon), \qquad \mathrm{F}^a_N(\varepsilon) := \frac{1}{N}\log \mathcal{Z}^a_{\varepsilon,N}, \tag{1.8}$$

where the existence of this limit (that will follow as a by-product of our approach) can be proven with a standard superadditivity argument. The basic observation is that the free energy is nonnegative. In fact, setting $\Omega^{\mathrm{w}} := [0,\infty)$ and $\Omega^{\mathrm{p}} := \mathbb{R}$, we have $\forall N \in \mathbb{N}$

$$\begin{aligned} \mathcal{Z}^a_{\varepsilon,N} &= \int \exp(-\mathcal{H}_{[-1,N+1]}(\varphi)) \prod_{i=1}^{N-1} (\varepsilon\delta_0(d\varphi_i) + d\varphi_i \mathbf{1}_{(\varphi_i\in\Omega^a)}) \\ &\geq \int \exp(-\mathcal{H}_{[-1,N+1]}(\varphi)) \prod_{i=1}^{N-1} d\varphi_i \mathbf{1}_{(\varphi_i\in\Omega^a)} = \mathcal{Z}^a_{0,N} \geq \frac{c_1}{N^{c_2}}, \end{aligned} \tag{1.9}$$

where $c_1, c_2$ are positive constants and the polynomial bound for $\mathcal{Z}^a_{0,N}$ (analogous to what happens in the gradient case; cf. [10]) is proven in (2.14). Therefore $\mathrm{F}^a(\varepsilon) \geq \mathrm{F}^a(0) = 0$ for every $\varepsilon \geq 0$. Since this lower bound has been obtained by ignoring the contribution of the paths that touch the defect line, one is led to the following

DEFINITION 1.1. For $a \in \{\mathrm{p}, \mathrm{w}\}$, the $a$-model $\{\mathbb{P}^a_{\varepsilon,N}\}_N$ is said to be localized if $\mathrm{F}^a(\varepsilon) > 0$.

The first problem is to understand for which values of $\varepsilon$ (if any) there is localization. Some considerations can be drawn easily. We introduce for convenience for $t \in \mathbb{R}$

$$\widetilde{\mathrm{F}}^a_N(t) := \mathrm{F}^a_N(e^t), \qquad \widetilde{\mathrm{F}}^a(t) := \mathrm{F}^a(e^t). \tag{1.10}$$

It is easy to show (see Appendix A) that $\widetilde{\mathrm{F}}_N^a(\cdot)$ is convex, therefore also $\widetilde{\mathrm{F}}^a(\cdot)$ is convex. In particular, the free energy $\mathrm{F}^a(\varepsilon) = \widetilde{\mathrm{F}}^a(\log \varepsilon)$ is a *continuous function*, as long as it is finite. $\mathrm{F}^a(\cdot)$ is also *nondecreasing*, because $\mathcal{Z}_{\varepsilon,N}^a$ is increasing in $\varepsilon$ [cf. the first line of (1.9)]. This observation implies that, for both $a \in \{\mathrm{p}, \mathrm{w}\}$, there is a critical value $\varepsilon_c^a \in [0, \infty]$ such that the $a$-model is localized if and only if $\varepsilon > \varepsilon_c^a$. Moreover $\varepsilon_c^{\mathrm{p}} \le \varepsilon_c^{\mathrm{w}}$, since $\mathcal{Z}_{\varepsilon,N}^{\mathrm{p}} \ge \mathcal{Z}_{\varepsilon,N}^{\mathrm{w}}$.

However, it is still not clear that a phase transition really exists, that is, that $\varepsilon_c^a \in (0, \infty)$. Indeed, in the gradient case the transition is nontrivial only for the wetting model, that is, $0 < \varepsilon_c^{\mathrm{w},\nabla} < \infty$ while $\varepsilon_c^{\mathrm{p},\nabla} = 0$; cf. [10, 13]. Our first theorem shows that in the Laplacian case both the pinning and the wetting models undergo a nontrivial transition, and gives further properties of the free energy $\mathrm{F}^a(\cdot)$.

THEOREM 1.2 (Localization transition). *The following relations hold:*

$$\varepsilon_c^{\mathrm{p}} \in (0, \infty), \qquad \varepsilon_c^{\mathrm{w}} \in (0, \infty), \qquad \varepsilon_c^{\mathrm{p}} < \varepsilon_c^{\mathrm{w}}.$$

*We have* $\mathrm{F}^a(\varepsilon) = 0$ *for* $\varepsilon \in [0, \varepsilon_c^a]$, *while* $0 < \mathrm{F}^a(\varepsilon) < \infty$ *for* $\varepsilon \in (\varepsilon_c^a, \infty)$, *and as* $\varepsilon \to \infty$

$$\mathrm{F}^a(\varepsilon) = \log \varepsilon (1 + o(1)), \qquad a \in \{\mathrm{p}, \mathrm{w}\}. \tag{1.11}$$

*Moreover the function* $\mathrm{F}^a(\varepsilon)$ *is real analytic on* $(\varepsilon_c^a, \infty)$.

One may ask why in the Laplacian case we have $\varepsilon_c^{\mathrm{p}} > 0$, unlike in the gradient case. Heuristically, we could say that the Laplacian interaction (1.3) describes a stiffer chain, more rigid to bending with respect to the gradient interaction, and therefore Laplacian models require a stronger reward in order to localize. Note in fact that in the Gaussian case $V(x) = x^2/(2\sigma^2) + \log \sqrt{2\pi\sigma^2}$ the ground state of the gradient interaction is just the horizontally flat line, whereas the Laplacian interaction favors rather *affine configurations*, penalizing curvature and bendings.

It is worth stressing that the free energy has a direct translation in terms of some path properties of the field. Defining the contact number $\ell_N$ by

$$\ell_N := \#\{i \in \{1, \ldots, N\} : \varphi_i = 0\}, \tag{1.12}$$

a simple computation (see Appendix A) shows that for every $\varepsilon > 0$ and $N \in \mathbb{N}$

$$\mathrm{D}_N^a(\varepsilon) := \mathbb{E}_{\varepsilon,N}^a\left(\frac{\ell_N}{N}\right) = (\widetilde{\mathrm{F}}_N^a)'(\log \varepsilon) = \varepsilon \cdot (\mathrm{F}_N^a)'(\varepsilon). \tag{1.13}$$

Then, introducing the nonrandom quantity $\mathrm{D}^a(\varepsilon) := \varepsilon \cdot (\mathrm{F}^a)'(\varepsilon)$ (which is well defined for $\varepsilon \ne \varepsilon_c^a$ by Theorem 1.2), a simple convexity argument shows that $\mathrm{D}_N^a(\varepsilon) \to \mathrm{D}^a(\varepsilon)$ as $N \to \infty$, for every $\varepsilon \ne \varepsilon_c^a$. Indeed much more can be said (see Appendix A):

- When $\varepsilon > \varepsilon_c^a$ we have that $\mathrm{D}^a(\varepsilon) > 0$, and for every $\delta > 0$ and $N \in \mathbb{N}$

$$\mathbb{P}^a_{\varepsilon,N}\left( \left| \frac{\ell_N}{N} - \mathrm{D}^a(\varepsilon) \right| > \delta \right) \le \exp(-c_3\, N), \tag{1.14}$$

  where $c_3$ is a positive constant. This shows that, when the $a$-model is localized according to Definition 1.1, its typical paths touch the defect line a positive fraction of times, equal to $\mathrm{D}^a(\varepsilon)$. Notice that, by (1.11) and convexity arguments, $\mathrm{D}^a(\varepsilon)$ converges to 1 as $\varepsilon \to \infty$, that is, a strong reward pins the field at the defect line in a very effective way (observe that $\ell_N/N \le 1$).

- On the other hand, when $\varepsilon < \varepsilon_c^a$ we have $\mathrm{D}^a(\varepsilon) = 0$ and for every $\delta > 0$ and $N \in \mathbb{N}$

$$\mathbb{P}^a_{\varepsilon,N}\left( \frac{\ell_N}{N} > \delta \right) \le \exp(-c_4\, N), \tag{1.15}$$

  where $c_4$ is a positive constant. Thus for $\varepsilon < \varepsilon_c^a$ the typical paths of the $a$-model touch the defect line only $o(N)$ times; when this happens it is customary to say that the model is *delocalized*.

What is left out from this analysis is the critical regime $\varepsilon = \varepsilon_c^a$. The behavior of the model in this case is sharply linked to the way in which the free energy $\mathrm{F}^a(\varepsilon)$ vanishes as $\varepsilon \downarrow \varepsilon_c^a$. If $\mathrm{F}^a(\cdot)$ is differentiable also at $\varepsilon = \varepsilon_c^a$ (transition of second or higher order), then $(\mathrm{F}^a)'(\varepsilon_c^a) = 0$ and relation 1.15 holds, that is, the $a$-model for $\varepsilon = \varepsilon_c^a$ is delocalized. The other possibility is that $\mathrm{F}^a(\cdot)$ is not differentiable at $\varepsilon = \varepsilon_c^a$ (transition of first order), which happens when the right-derivative is positive: $(\mathrm{F}^a)'_+(\varepsilon_c^a) > 0$. In this case the behavior of $\mathbb{P}^a_{\varepsilon,N}$ for large $N$ may depend on the choice of the boundary conditions.

We first consider the critical regime for the wetting model, where the transition turns out to be of *first order*. Recall the definition 1.13 of $\mathrm{D}^a_N(\varepsilon)$.

THEOREM 1.3 (Critical wetting model). *For the wetting model we have*

$$\liminf_{N\to\infty} \mathrm{D}^{\mathrm{w}}_N(\varepsilon_c^{\mathrm{w}}) > 0. \tag{1.16}$$

*Therefore* $(\mathrm{F}^{\mathrm{w}})'_+(\varepsilon_c^{\mathrm{w}}) > 0$ *and the phase transition is of first order.*

Notice that (1.13) and (1.16) yield

$$\liminf_{N\to\infty} \mathbb{E}^{\mathrm{w}}_{\varepsilon_c^{\mathrm{w}},N}\left( \frac{\ell_N}{N} \right) > 0,$$

and in this sense the wetting model at the critical point exhibits a *localized* behavior. This is in sharp contrast with the gradient case, where it is well known that the wetting model at criticality is delocalized and in fact the transition is of *second order*; cf. [9, 10, 11, 18]. The emergence of a first-order transition in the case of Laplacian interaction is interesting in view of the possible applications of $\mathbb{P}^{\mathrm{w}}_{\varepsilon,N}$ as a model for the *DNA denaturation transition*, where the nonnegative field $\{\varphi_i\}_i$ describes the distance between the two DNA strands. In fact for the DNA denaturation something close to a first-order phase transition is experimentally observed; we refer to [13], Section 1.4, for a detailed discussion (cf. also [19, 27]).

Finally we consider the critical pinning model, where the transition is of second order.

THEOREM 1.4 (Critical pinning model). *For the pinning model we have*

$$\limsup_{\varepsilon\downarrow\varepsilon_c^{\mathrm{p}}}\limsup_{N\to\infty} \mathrm{D}^{\mathrm{p}}_N(\varepsilon)=0. \tag{1.17}$$

*Therefore* $\mathrm{F}^{\mathrm{p}}(\varepsilon)$ *is differentiable at* $\varepsilon=\varepsilon_c^{\mathrm{p}}$ *and* $(\mathrm{F}^{\mathrm{p}})'(\varepsilon_c^{\mathrm{p}})=0$. *Moreover there exists* $c_5>0$ *such that for* $\delta$ *sufficiently small we have*

$$\mathrm{F}^{\mathrm{p}}(\varepsilon_c^{\mathrm{p}}+\delta,0)\geq c_5\frac{\delta}{\log 1/\delta}, \tag{1.18}$$

*that is, the transition is exactly of second order.*

Although the relation $(\mathrm{F}^{\mathrm{p}})'(\varepsilon_c^{\mathrm{p}})=0$ yields $\ell_N=o(N)$, in a delocalized fashion, the pinning model at $\varepsilon=\varepsilon_c^{\mathrm{p}}$ is actually somewhat borderline between localization and delocalization, as (1.18) suggests and as we point out in the next paragraph.

1.3. *Further path results.* A direct application of the techniques that we develop in this paper yields further path properties of the field. Let us introduce the *maximal gap*

$$\Delta_N:=\max\{n\leq N:\varphi_{k+1}\neq 0,\varphi_{k+2}\neq 0,\dots,\varphi_{k+n}\neq 0 \text{ for some } k\leq N-n\}.$$

One can show that, for both $a\in\{\mathrm{p},\mathrm{w}\}$ and for $\varepsilon>\varepsilon_c^a$, the following relations hold:

$$\begin{gathered}\forall\delta>0\colon \lim_{N\to\infty}\mathbb{P}^a_{\varepsilon,N}\left(\frac{\Delta_N}{N}\geq\delta\right)=0,\\ \lim_{L\to\infty}\limsup_{N\to\infty}\max_{i=1,\dots,N-1}\mathbb{P}^a_{\varepsilon,N}(|\varphi_i|\geq L)=0.\end{gathered} \tag{1.19}$$

In particular for $\varepsilon>\varepsilon_c^a$ each component $\varphi_i$ of the field is at finite distance from the defect line and this is a clear localization path statement. On the other hand, in the pinning case $a=\mathrm{p}$ we can strengthen (1.15) to the following relation: for every $\varepsilon<\varepsilon_c^{\mathrm{p}}$

$$\lim_{L\to\infty}\limsup_{N\to\infty}\mathbb{P}^{\mathrm{p}}_{\varepsilon,N}(\Delta_N\leq N-L)=0, \tag{1.20}$$

that is, for $\varepsilon < \varepsilon_c^{\mathrm{p}}$ the field touches the defect line at a finite number of sites, all at *finite distance* from the boundary points $\{0, N\}$. We expect that the same relation holds true also in the wetting case $a = \mathrm{w}$, but at present we cannot prove it: what is missing are more precise estimates on the entropic repulsion problem; see Section 1.5 for a detailed discussion. It is interesting to note that we can prove that the first relation in (1.19) holds true also in the pinning case $a = \mathrm{p}$ at the critical point $\varepsilon = \varepsilon_c^{\mathrm{p}}$, and this shows that the pinning model at criticality has also features of localized behavior.

We do not give an explicit proof of the above relations in this paper, both for conciseness and because in a second paper [7] we focus on the scaling limits of the pinning model, obtaining (de)localization path statements that are much more precise than (1.19) and (1.20) [under stronger assumptions on the potential $V(\cdot)$]. We show in particular that for all $\varepsilon \in (0, \varepsilon_c^{\mathrm{p}})$ the natural rescaling of $\mathbb{P}_{\varepsilon,N}^{\mathrm{p}}$ converges in distribution in $C([0,1])$ to the same limit that one obtains in the free case $\varepsilon = 0$, that is, the integral process of a Brownian bridge. On the other hand, for every $\varepsilon \geq \varepsilon_c^{\mathrm{p}}$ the natural rescaling of $\mathbb{P}_{\varepsilon,N}^{\mathrm{p}}$ yields the trivial process which is identically zero. We stress that $\varepsilon = \varepsilon_c^{\mathrm{p}}$ is included in the last statement; this is in sharp contrast with the gradient case, where both the pinning and the wetting models at criticality have a nontrivial scaling limit, respectively the Brownian bridge and the reflected Brownian bridge. This shows again the peculiarity of the critical pinning model in the Laplacian case. Indeed, by lowering the scaling constants with suitable logarithmic corrections, we are able to extract a nontrivial scaling limit (in a distributional sense) for the law $\mathbb{P}_{\varepsilon_c^{\mathrm{p}},N}^{\mathrm{p}}$ in terms of a symmetric stable Lévy process. We stress that the techniques and results of the present paper play a crucial role for [7].

1.4. *Outline of the paper: approach and techniques.* Although our main results are about the free energy, the core of our approach is a precise pathwise description of the field based on Markov renewal theory. In analogy to [9, 10] and especially to [8], we would like to stress the power of (Markov) renewal theory techniques for the study of $(1+1)$-dimensional linear chain models. The other basic techniques that we use are local limit theorems and an infinite-dimensional version of the Perron–Frobenius theorem. Let us describe more in detail the structure of the paper.

In Section 2 we study the pinning and wetting models in the free case $\varepsilon = 0$, showing that these models are sharply linked to the integral of a random walk. More precisely, let $\{Y_n\}_{n \geq 0}$ denote a random walk starting at zero and with step law $\mathbf{P}(Y_1 \in dx) = \exp(-V(x))\, dx$ [the walk has zero mean and finite variance by (1.2)] and let us denote by $Z_n := Y_1 + \cdots + Y_n$ the corresponding *integrated random walk process*. In Proposition 2.2 we show that the law $\mathbb{P}_{0,N}^{a}$ is nothing but a bridge of length $N$ of the process $\{Z_n\}_n$, with the further conditioning to stay nonnegative in the wetting case $a = \mathrm{w}$. Therefore we focus on the asymptotic properties of the process $\{Z_n\}_n$, obtaining a basic local limit theorem (cf. Proposition 2.3) and some bounds for the probability that $\{Z_n\}_n$ stays positive (connected to the problem of entropic repulsion that we discuss below; cf. Section 1.5).

In Section 3, which is the core of the paper, we show that for $\varepsilon > 0$ the law $\mathbb{P}^a_{\varepsilon,N}$ admits a crucial description in terms of Markov renewal theory. More precisely, we show that the zeros of the field $\{i \le N : \varphi_i = 0\}$ under $\mathbb{P}^a_{\varepsilon,N}$ are distributed according to the law of a (hidden) Markov renewal process conditioned to hit $\{N, N+1\}$; cf. Proposition 3.1. We thus obtain an explicit expression for the partition function $\mathcal{Z}^a_{\varepsilon,N}$ in terms of this Markov renewal process, which is the key to our main results.

Section 4 is devoted to proving some analytical results that underlie the construction of the Markov renewal process appearing in Section 3. The main tool is an infinite-dimensional version of the classical Perron–Frobenius theorem (cf. [33]), and a basic role is played by the asymptotic estimates obtained in Section 2. A by-product of this analysis is an explicit formula [cf. 4.10], that links $\mathrm{F}^a(\cdot)$ and $\varepsilon^a_c$ to the spectral radius of a suitable integral operator and that will be exploited later.

Sections 5, 6 and 7 contain the proofs of Theorems 1.2, 1.3 and 1.4, respectively. In view of the description given in Section 3, all the results to prove can be rephrased in the language of Markov renewal theory. The proofs are then carried out exploiting the asymptotic estimates derived in Sections 2 and 4 together with some algebraic manipulation of the kernel that gives the law of the hidden Markov renewal process. Finally, the Appendixes contain the proof of some technical results.

We conclude the introduction by discussing briefly two interesting problems that are linked to our models, namely the *entropic repulsion* in Section 1.5 and the *smoothing effect of disorder* in Section 1.6. Finally, Section 1.7 contains some recurrent notations, especially about kernels, used throughout the paper.

1.5. *Entropic repulsion.* We recall that $(\{Y_n\}_{n\ge 0}, \mathbf{P})$ is the random walk with step $\mathbf{P}(Y_1 \in dx) = e^{-V(x)}\,dx$ and that $Z_n = Y_1 + \cdots + Y_n$. The analysis of the wetting model requires estimating the decay as $N \to \infty$ of the probabilities $\mathbf{P}(\Omega_N^+)$ and $\mathbf{P}(\Omega_N^+ \mid Z_{N+1} = 0, Z_{N+2} = 0)$, where we set $\Omega_N^+ := \{Z_1 \ge 0, \ldots, Z_N \ge 0\}$. This type of problem is known in the literature as *entropic repulsion* and it has received a lot of attention; see [32] for a recent overview. In the Laplacian case that we consider here, this problem has been solved in the Gaussian setting (i.e., when $V(x) = x^2/(2\sigma^2) + \log\sqrt{2\pi\sigma^2}$) in $(d+1)$-dimension with $d \ge 5$; cf. [22, 28]. Little is known in the $(1+1)$-dimensional setting, apart from the following result of Sinai's [29] in the very special case when $\{Y_n\}_n$ is the simple random walk on $\mathbb{Z}$:

$$\frac{c}{N^{1/4}} \leq \mathbf{P}(\Omega_N^+) \leq \frac{C}{N^{1/4}}, \tag{1.21}$$

where $c, C$ are positive constants. The proof of this bound relies on the exact combinatorial results available in the simple random walk case and it appears difficult to extend it to our situation. We point out that the same exponent $1/4$ appears in related continuous models dealing with the integral of Brownian motion; cf. [23, 25]. Based on Sinai's result, which we believe to hold for general random walks with zero mean and finite variance, we expect that for the bridge case one should have the bound

$$\frac{c}{N^{1/2}} \leq \mathbf{P}(\Omega_N^+ \mid Z_{N+1} = 0, Z_{N+2} = 0) \leq \frac{C}{N^{1/2}}. \tag{1.22}$$

We cannot derive precise bounds as (1.21) and (1.22); however, for the purpose of this paper the following weaker result suffices:

PROPOSITION 1.5. *There exist constants $c, C, c_- > 0$ and $c_+ > 1$ such that $\forall N \in \mathbb{N}$*

$$\frac{c}{N^{c_-}} \leq \mathbf{P}(\Omega_N^+) \leq \frac{C}{(\log N)^{c_+}}, \tag{1.23}$$

$$\frac{c}{N^{c_-}} \leq \mathbf{P}(\Omega_N^+ \mid Z_{N+1} = 0, Z_{N+2} = 0) \leq \frac{C}{(\log N)^{c_+}}. \tag{1.24}$$

We prove this proposition in Appendix C. We stress that the most delicate point is the upper bound in (1.23): the idea of the proof is to dilute the system on superexponentially spaced times. We also point out that in the Gaussian case, that is, when $V(x) = x^2/(2\sigma^2) + \log\sqrt{2\pi\sigma^2}$, the upper bound can be easily strengthened to $\mathbf{P}(\Omega_N^+) \leq (const.)/N^a$, for some $a > 0$, by diluting the system on exponentially spaced times and using the results in the paper [24].

1.6. *Smoothing effect of quenched disorder.* The models we consider in this work are *homogeneous*, in that the reward $\varepsilon$ is deterministic and constant. However, one can define in a natural way a disordered version of our model, where the reward is itself random and may vary from site to site. We stress that disordered pinning models have attracted a lot of attention recently; cf. [1, 2, 16, 30, 31] (see also [13] for an overview).

Let us be more precise: we take a sequence $\omega = \{\omega_n\}_{n\in\mathbb{N}}$ of i.i.d. standard Gaussian random variables, defined on some probability space $(\mathcal{S}, P)$, and we denote by $\mathrm{M}(\beta) := E(\exp(\beta\omega_1)) = \exp(\beta^2/2)$ the corresponding moment generating function (we focus on Gaussian variables only for the sake of simplicity; we could more generally take variables with zero mean, unit variance and finite exponential moments). Then, given a realization $\omega$, we introduce for $\beta \geq 0$ a *quenched* version $\mathbb{P}^a_{\varepsilon,\beta,\omega,N}$ of our polymer measure, just replacing the reward $\varepsilon$ by $\varepsilon \exp(\beta\omega_i)$:

$$\mathbb{P}^a_{\varepsilon,\beta,\omega,N}(d\varphi_1 \cdots d\varphi_{N-1}) := \frac{\exp(-\mathcal{H}_{[-1,N+1]}(\varphi))}{\mathcal{Z}^a_{\varepsilon,\beta,\omega,N}} \prod_{i=1}^{N-1} (\varepsilon e^{\beta\omega_i} \delta_0(d\varphi_i) + d\varphi_i \mathbf{1}_{\{\varphi_i \in \Omega^a\}}),$$

where of course $a \in \{\mathrm{p}, \mathrm{w}\}$ and $\Omega^{\mathrm{p}} = \mathbb{R}$ while $\Omega^{\mathrm{w}} = \mathbb{R}^+$.

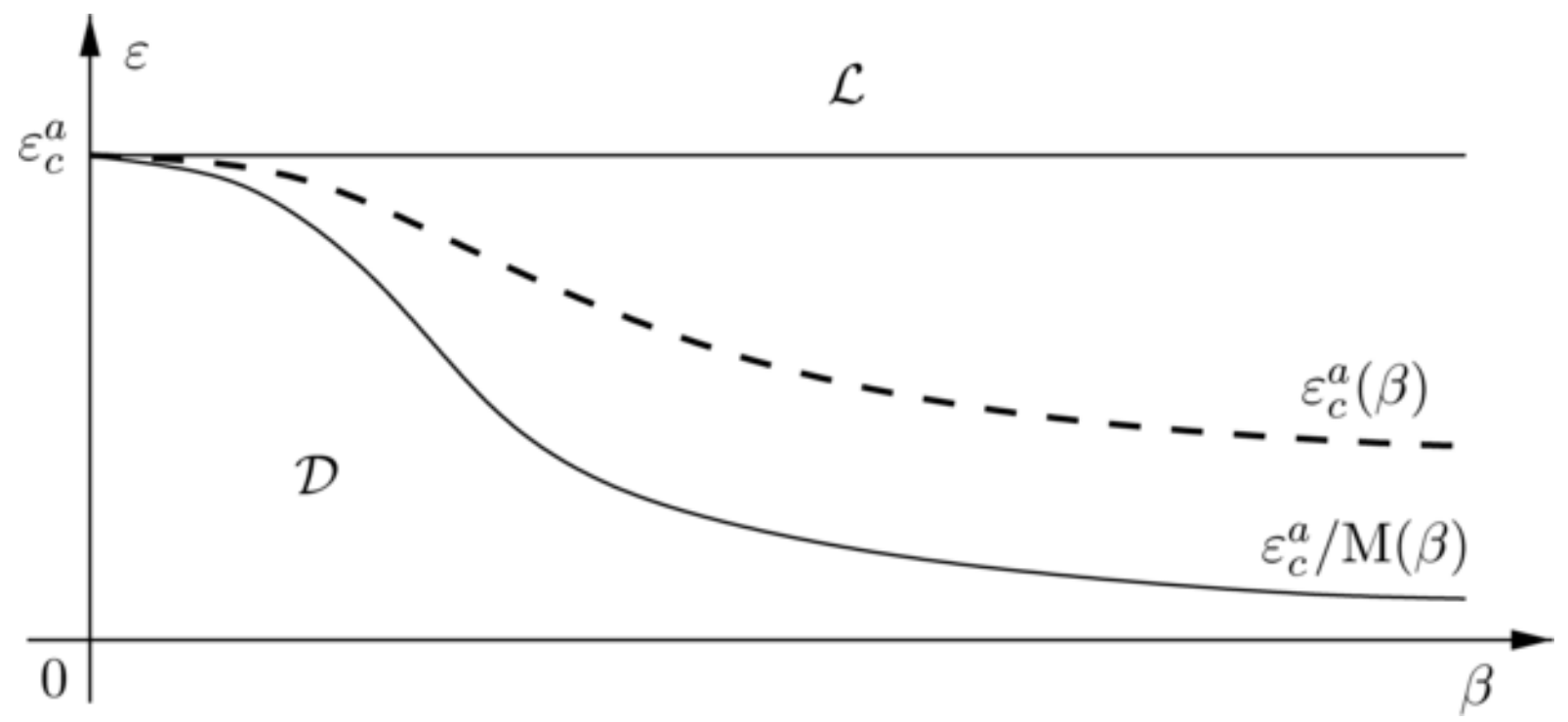

FIG. 2. *The critical value* $\varepsilon_c^a(\beta)$ *as a function of* $\beta$, *with the upper and lower bound given in (1.25).* $\mathcal{L} = \{(\beta,\varepsilon) : \varepsilon > \varepsilon_c^a(\beta)\}$ *and* $\mathcal{D} = \{(\beta,\varepsilon) : \varepsilon < \varepsilon_c^a(\beta)\}$ *represent the localized and delocalized region, respectively.*

The analysis of the model goes along the same line as in the homogeneous case. Namely, the free energy $\mathrm{F}^a(\varepsilon,\beta)$ is defined by

$$\mathrm{F}^a(\varepsilon,\beta) := \lim_{N\to\infty} \mathrm{F}^a_N(\varepsilon,\beta,\omega), \qquad \mathrm{F}^a_N(\varepsilon,\beta,\omega) := \frac{1}{N} \log \mathcal{Z}^a_{\varepsilon,\beta,\omega,N},$$

where the existence of this limit, both $P(d\omega)$-a.s. and in $L_1(dP)$, follows by Kingman's Superadditive Ergodic Theorem (cf. [13], Section 4.4.2), and moreover $\mathrm{F}^a(\varepsilon,\beta)$ does not depend on $\omega$, that is, it is $P(d\omega)$-a.s. a constant (*self-averaging* of the free energy).

Restricting to the paths that do not touch the defect line, exactly like in (1.9), yields the lower bound $\mathrm{F}^a(\varepsilon,\beta) \geq 0$ for all $\varepsilon,\beta$. Therefore, as in the homogeneous case, we say that the $a$-model with parameters $(\varepsilon,\beta)$ is *localized* if $\mathrm{F}^a(\varepsilon,\beta) > 0$. By the monotonicity in $\varepsilon$ it follows that, for every value of $\beta \geq 0$, there is a critical value $\varepsilon_c^a(\beta) \in [0,\infty]$ such that the $a$-model is localized if and only if $\varepsilon > \varepsilon_c^a(\beta)$. Of course $\varepsilon_c^a(0)$ equals the value $\varepsilon_c^a$, defined before Theorem 1.2. Using convexity arguments and Jensen's inequality (*annealed bound*), in close analogy with [13], Proposition 5.1, one proves easily that

$$\mathrm{F}^a(\varepsilon,0) \leq \mathrm{F}^a(\varepsilon,\beta) \leq \mathrm{F}^a(\varepsilon\mathrm{M}(\beta),0) \qquad \forall \beta \geq 0.$$

This yields immediately that for every $\beta \geq 0$

$$\frac{\varepsilon_c^a}{\mathrm{M}(\beta)} \leq \varepsilon_c^a(\beta) \leq \varepsilon_c^a; \tag{1.25}$$

see Figure 2. In particular, $0 < \varepsilon_c^a(\beta) < \infty$ for every $\beta \geq 0$, by Theorem 1.2. By convexity arguments one also shows that the curve $\beta \mapsto \varepsilon_c^a(\beta)$ is continuous and nonincreasing.

The picture emerging from these simple considerations is that, also in presence of disorder, there is still a nontrivial phase transition: by increasing the parameter $\varepsilon$ one passes from the delocalized region $\varepsilon < \varepsilon_c^a(\beta)$ to the localized region $\varepsilon > \varepsilon_c^a(\beta)$. When disorder is absent, that is, for $\beta = 0$, the behavior of the free energy near the phase transition critical point $\varepsilon_c^a = \varepsilon_c^a(0)$ is described by Theorems 1.3 and 1.4, for the wetting and the pinning models, respectively, and it is natural to ask what happens when disorder is switched on, that is, for $\beta > 0$. It turns out that the presence of disorder has a strong smoothing effect on the phase transition. More precisely, one can prove the following result.

PROPOSITION 1.6. *For every $\beta > 0$ there is a constant $C_\beta \in (0, \infty)$ such that for both $a = \mathrm{p}, \mathrm{w}$ the following inequality holds true for small $\delta > 0$:*

$$\mathrm{F}^a(\varepsilon_c^a(\beta) + \delta, \beta) \leq C_\beta \delta^2. \tag{1.26}$$

This means that, as soon as disorder is present, the phase transition is at least of second order, with *bounded* second derivative of the free energy. Let us stress the sharp difference with the results we have found for the homogeneous case $\beta = 0$:

- in the wetting case, by Theorem 1.3 the transition is of first order and in fact

$$\mathrm{F}^{\mathrm{w}}(\varepsilon_c^{\mathrm{w}} + \delta, 0) \geq c_6 \delta,$$

  for small $\delta > 0$ and for some positive constant $c_6$, in contrast with (1.26);
- even in the pinning case, where the transition is of second order, (1.18) shows that the second derivative of the free energy from the right at $\varepsilon = \varepsilon_c^{\mathrm{p}}$ is *infinite*, in contrast with (1.26).

We omit a detailed proof of Proposition 1.6, because it follows exactly the same line as that of [15], Theorem 2.1, and we just sketch it. The idea is to search for *long*, *atypical stretches* in the disorder sequence $\omega$ and to give an explicit *lower bound* on the partition function, by making the field come back to the defect line only in correspondence to these stretches. Then, in a clever way, one extracts from this lower bound an *upper bound* on the free energy near the critical point (for a beautiful heuristic description see [14]).

The original proof of [15], Theorem 2.1, holds for a general class of disordered pinning models, for which the free models (i.e., without pinning) enjoy the following properties:

1. the epochs of the returns to zero have a *genuine renewal structure*;
2. the corresponding interarrival time probability has a power-law lower bound.

This class of models includes in particular $(1+1)$-dimensional pinning and wetting models with *gradient interaction*. As we already mentioned, in our Laplacian-interaction case the returns to zero have only a *Markov* renewal structure. Nevertheless, the set $\{i \in \mathbb{N} : \varphi_i = \varphi_{i-1} = 0\}$ of *double zeros* has a genuine renewal structure, as we show in Section 5.1, and this is sufficient for the proof, because making the field come back to the defect line with double zeros gives a lower bound. Finally, the power-law lower bound for the interarrival time probability in the free model (i.e., when $\varepsilon = 0$) is given by Proposition 2.3 in the pinning case and by Proposition 1.5 in the wetting case. This shows again the usefulness of bounds like (1.23) and (1.24), despite the fact that they are not sharp.

1.7. *Some recurrent notation.* Throughout the paper, generic positive and finite constants will be denoted by $(const.)$, $(const.')$. For us $\mathbb{N} = \{1, 2, \ldots\}$, $\mathbb{Z}^+ = \mathbb{N} \cup \{0\}$ and $\mathbb{R}^+ := [0, \infty)$. Given two positive sequences $(a_n)$, $(b_n)$, by $a_n \sim b_n$ we mean that $a_n / b_n \to 1$ as $n \to \infty$. For $x \in \mathbb{R}$ we denote as usual by $\lfloor x \rfloor := \max\{n \in \mathbb{Z} : n \le x\}$ its integer part.

In this paper we deal with kernels of two kinds. Kernels of the first kind are just $\sigma$-finite kernels on $\mathbb{R}$, that is, functions $A_{\cdot,\cdot} : \mathbb{R} \times \mathcal{B}(\mathbb{R}) \to \mathbb{R}^+$, where $\mathcal{B}(\mathbb{R})$ denotes the Borel $\sigma$-field of $\mathbb{R}$, such that $A_{x,\cdot}$ is a $\sigma$-finite Borel measure on $\mathbb{R}$ for every $x \in \mathbb{R}$ and $A_{\cdot,F}$ is a Borel function for every $F \in \mathcal{B}(\mathbb{R})$. Given two such kernels $A_{x,dy}$, $B_{x,dy}$, their *composition* is denoted as usual by $(A \circ B)_{x,dy} := \int_{z \in \mathbb{R}} A_{x,dz} B_{z,dy}$ and $A^{\circ k}_{x,dy}$ denotes the $k$-fold composition of $A$ with itself, where $A^{\circ 0}_{x,dy} := \delta_x(dy)$. We also use the standard notation

$$(1 - A)^{-1}_{x,dy} := \sum_{k=0}^{\infty} A^{\circ k}_{x,dy},$$

which of course in general may be infinite.

The second kind of kernel is obtained by letting a kernel of the first kind depend on the further parameter $n \in \mathbb{Z}^+$, that is, we consider objects of the form $\mathsf{A}_{x,dy}(n)$ with $x, y \in \mathbb{R}$ and $n \in \mathbb{Z}^+$. Given two such kernels $\mathsf{A}_{x,dy}(n)$, $\mathsf{B}_{x,dy}(n)$, we define their *convolution* by

$$(\mathsf{A} * \mathsf{B})_{x,dy}(n) := \sum_{m=0}^{n} (\mathsf{A}(m) \circ \mathsf{B}(n-m))_{x,dy}$$

$$= \sum_{m=0}^{n} \int_{z\in\mathbb{R}} \mathsf{A}_{x,dz}(m)\mathsf{B}_{z,dy}(n-m),$$

and the $k$-fold convolution of the kernel $\mathsf{A}_{x,dy}(n)$ with itself will be denoted by $\mathsf{A}^{*k}_{x,dy}(n)$, where by definition $\mathsf{A}^{*0}_{x,dy}(n) := \delta_0(dy)\mathbf{1}_{(n=0)}$. Finally, given two kernels $\mathsf{A}_{x,dy}(n)$ and $B_{x,dy}$ and a positive sequence $(a_n)$, we will write

$$\mathsf{A}_{x,dy}(n) \sim \frac{B_{x,dy}}{a_n} \qquad (n\to\infty)$$

to mean $\mathsf{A}_{x,F}(n) \sim B_{x,F}/a_n$ as $n\to\infty$, $\forall x\in\mathbb{R}$ and for every *bounded* Borel set $F\subset\mathbb{R}$.

**2. The free case $\varepsilon = 0$: a random walk viewpoint.** In this section we study in detail the free laws $\mathbb{P}^{\mathrm{p}}_{0,N}$ and $\mathbb{P}^{\mathrm{w}}_{0,N}$ and their link with the integral of a random walk. The main results are a basic local limit theorem and some asymptotic estimates.

2.1. *Integrated random walk.* Given $a,b\in\mathbb{R}$, let $(\Omega,\mathcal{F},\mathbf{P}=\mathbf{P}^{(a,b)})$ be a probability space on which are defined the processes $\{X_i\}_{i\in\mathbb{N}}$, $\{Y_i\}_{i\in\mathbb{Z}^+}$ and $\{Z_i\}_{i\in\mathbb{Z}^+}$ with the following properties:

- $\{X_i\}_{i\in\mathbb{N}}$ is a sequence of independent and identically distributed random variables, with marginal laws $X_1\sim\exp(-V(x))\,dx$. We recall that by our assumptions on $V(\cdot)$ it follows that $\mathbf{E}(X_1)=0$ and $\mathbf{E}(X_1^2)=\sigma^2\in(0,\infty)$; cf. (1.2).
- $\{Y_i\}_{i\in\mathbb{Z}^+}$ is the random walk associated to $\{X_i\}$, with starting point $a$, that is,

$$Y_0 = a, \qquad Y_n = a + X_1 + \cdots + X_n. \tag{2.1}$$

- $\{Z_i\}_{i\in\mathbb{Z}^+}$ is the *integrated random walk process* with initial value $b$; that is, $Z_0=b$ and for $n\in\mathbb{N}$

$$Z_n = b + Y_1 + \cdots + Y_n = b + na + nX_1 + (n-1)X_2 + \cdots + X_n. \tag{2.2}$$

From (2.1) and (2.2) it follows that

$$\{(Y_n,Z_n)\}_n \text{ under } \mathbf{P}^{(a,b)} \overset{d}{=} \{(Y_n+a, Z_n+b+na)\}_n \text{ under } \mathbf{P}^{(0,0)}. \tag{2.3}$$

The marginal distributions of the process $\{Z_n\}_n$ are specified in the following lemma.

Lemma 2.1. *For every $n\in\mathbb{N}$, the law of the vector $(Z_1,\ldots,Z_n)$ under $\mathbf{P}^{(a,b)}$ is given by*

$$\begin{aligned} &\mathbf{P}^{(a,b)}((Z_1,\ldots,Z_n)\in(dz_1,\ldots,dz_n)) \\ &\quad = \exp(-\mathcal{H}_{[-1,n]}(z_{-1},z_0,z_1,\ldots,z_n))\prod_{i=1}^n dz_i, \end{aligned} \tag{2.4}$$

*where we set* $z_{-1} := b - a$ *and* $z_0 := b$.

PROOF. By (2.2) we have $Y_n = Z_n - Z_{n-1}$ for $n \geq 1$. Then, setting $y_i := z_i - z_{i-1}$ for $i \geq 2$ and $y_1 := z_1 - b$, it suffices to show that, under the measure given by the r.h.s. of (2.4), the variables $(y_i)_{i=1,\dots,n}$ are distributed like the first $n$ steps of a random walk starting at $a$ and with step law $\exp(-V(x))\,dx$. But for this it suffices to rewrite the Hamiltonian as

$$\begin{aligned}\mathcal{H}_{[-1,n]}(z) &= V((z_1 - b) - (b - (b - a))) + \sum_{i=1}^{n-1} V((z_{i+1} - z_i) - (z_i - z_{i-1}))\\ &= V(y_1 - a) + \sum_{i=1}^{n-1} V(y_{i+1} - y_i),\end{aligned}$$

and the proof is completed. □

By construction $\{(Y_n, Z_n)\}_{n\in\mathbb{Z}^+}$ under $\mathbf{P}^{(a,b)}$ is a Markov process with starting values $Y_0 = a$, $Z_0 = b$. On the other hand, the process $\{Z_n\}_n$ alone is not a Markov process; it is rather a process with *finite memory* $m = 2$, that is, for every $n \in \mathbb{N}$

$$\begin{aligned}\mathbf{P}^{(a,b)}(\{Z_{n+k}\}_{k\geq 0} \in \cdot \mid Z_i, i \leq n) &= \mathbf{P}^{(a,b)}(\{Z_{n+k}\}_{k\geq 0} \in \cdot \mid Z_{n-1}, Z_n)\\ &= \mathbf{P}^{(Z_n - Z_{n-1}, Z_n)}(\{Z_k\}_{k\geq 0} \in \cdot),\end{aligned} \tag{2.5}$$

as follows from Lemma 2.1. For this reason the law $\mathbf{P}^{(a,b)}$ may be viewed as

$$\mathbf{P}^{(a,b)} = \mathbf{P}(\cdot | Z_{-1} = b - a, Z_0 = b). \tag{2.6}$$

2.2. *The link with* $\mathbb{P}^a_{0,N}$. In the r.h.s. of (2.4) we see exactly the same density appearing in the definitions (1.5) and (1.7) of our models $\mathbb{P}^a_{0,N}$. As an immediate consequence we have the following proposition, which states that $\mathbb{P}^a_{0,N}$ is nothing but a bridge of the process $\{Z_n\}_n$ for $a = \mathrm{p}$, with the further constraint to stay nonnegative for $a = \mathrm{w}$.

PROPOSITION 2.2. *The following statements hold:*

(1) *The pinning model* $\mathbb{P}^{\mathrm{p}}_{0,N}$ *is the law of the vector* $(Z_1, \dots, Z_{N-1})$ *under the measure* $\mathbf{P}^{(0,0)}(\cdot \mid Z_N = 0, Z_{N+1} = 0)$. *The partition function* $\mathcal{Z}^{\mathrm{p}}_{0,N}$ *is the value at* $(0,0)$ *of the density of the vector* $(Y_{N+1}, Z_{N+1})$ *under the law* $\mathbf{P}^{(0,0)}$.

(2) *Setting* $\Omega^+_n := \{S_1 \geq 0, \dots, S_n \geq 0\}$, *the wetting model* $\mathbb{P}^{\mathrm{w}}_{0,N}$ *is the law of the vector* $(Z_1, \dots, Z_{N-1})$ *under the measure* $\mathbf{P}^{(0,0)}(\cdot | \Omega^+_{N-1}, Z_N = 0, Z_{N+1} = 0)$. *The partition function* $\mathcal{Z}^{\mathrm{w}}_{0,N}$ *is the value at* $(0,0)$ *of the density of the vector* $(Y_{N+1}, Z_{N+1})$ *under the law* $\mathbf{P}^{(0,0)}(\cdot \mid \Omega^+_{N-1})$, *multiplied by* $\mathbf{P}^{(0,0)}(\Omega^+_{N-1})$.

For the statement on the partition function, observe that the density at $(0,0)$ of the vector $(Y_{N+1}, Z_{N+1})$ coincides with the one of the vector $(Z_N, Z_{N+1})$, since $Y_{N+1} = Z_{N+1} - Z_N$.

2.3. *A local limit theorem.* In view of Proposition 2.2, we study the asymptotic behavior as $n \to \infty$ of the vector $(Y_n, Z_n)$ under the law $\mathbf{P}^{(a,b)}$.

Let us denote by $\{B_t\}_{t\in[0,1]}$ a standard Brownian motion and by $\{I_t\}_{t\in[0,1]}$ its integral process $I_t := \int_0^t B_s\, ds$. A simple application of Donsker's invariance principle shows that the vector $(Y_n/(\sigma\sqrt{n}), Z_n/(\sigma n^{3/2}))$ under $\mathbf{P}^{(0,0)}$ converges in distribution as $n \to \infty$ toward the law of the centered Gaussian vector $(B_1, I_1)$, whose density $g(y,z)$ is

$$g(y,z) = \frac{6}{\pi} \exp(-2y^2 - 6z^2 + 6yz). \tag{2.7}$$

We want to reinforce this convergence in the form of a local limit theorem. To this purpose, we introduce the density of $(Y_n, Z_n)$ under $\mathbf{P}^{(a,b)}$, setting for $n \geq 2$

$$\varphi_n^{(a,b)}(y,z) = \frac{\mathbf{P}^{(a,b)}(\,(Y_n, Z_n) \in (dy, dz))}{dy\, dz}. \tag{2.8}$$

From 2.3 it follows that

$$\varphi_n^{(a,b)}(y,z) = \varphi_n^{(0,0)}(y-a, z-b-na), \tag{2.9}$$

hence it suffices to focus on $\varphi_n^{(0,0)}(\cdot,\cdot)$. We set for short $\varphi_n^{(0,0)}(\mathbb{R}, z) := \int_{\mathbb{R}} \varphi_n^{(0,0)}(y,z)\, dy$, that is, the density of $Z_n$ under $\mathbf{P}^{(0,0)}$, and $g(\mathbb{R}, z) := \int_{\mathbb{R}} g(y,z)\, dy$. We are ready to state the main result of this section.

Proposition 2.3 (Local limit theorem). *The following relations hold as $n \to \infty$:*

$$\begin{aligned} &\sup_{(y,z)\in\mathbb{R}^2} |\sigma^2 n^2 \varphi_n^{(0,0)}(\sigma\sqrt{n}\, y, \sigma n^{3/2} z) - g(y,z)| \to 0, \\ &\sup_{z\in\mathbb{R}} |\sigma n^{3/2} \varphi_n^{(0,0)}(\mathbb{R}, \sigma n^{3/2}\, z) - g(\mathbb{R}, z)| \to 0. \end{aligned} \tag{2.10}$$

The proof, based on Fourier analysis, is deferred to Appendix B. We stress that this result retains a crucial importance in the rest of the paper. Notice that an analogous local limit theorem holds also for $\varphi_n^{(0,0)}(y, \mathbb{R})$, that is, the density of $Y_n$, but we do not state it explicitly because we will not need it.

2.4. *The positivity constraint.* To deal with the wetting model we need to study the law of the random vector $(Y_n, Z_n)$, or equivalently of $(Z_{n-1}, Z_n)$, conditionally on the event $\Omega^+_{n-2} = \{Z_1 \geq 0, \dots, Z_{n-2} \geq 0\}$. To this purpose we set for $x, y \in \mathbb{R}$ and $n \geq 3$

$$w_{x,y}(n) := \mathbf{1}_{(x\geq 0,\, y\geq 0)} \cdot \mathbf{P}^{(-x,0)}(\Omega^+_{n-2} \,|\, Z_{n-1} = y, Z_n = 0), \tag{2.11}$$

while for $n = 1, 2$ we simply set $w_{x,y}(n) := \mathbf{1}_{(x\geq 0, y\geq 0)}$. We are interested in the rate of decay of $w_{x,y}(n)$ as $n \to \infty$. To this purpose we claim that there exists a positive constant $c_+ > 1$ such that the following upper bound holds: for all $n \in \mathbb{N}$

$$\sup_{0\leq x,y\leq \sqrt{n}} w_{x,y}(n) \leq \frac{(const.)}{(\log n)^{c_+}}. \tag{2.12}$$

Moreover we have the following lower bound for $x, y = 0$ and for all $n \in \mathbb{N}$:

$$w_{0,0}(n) \geq \frac{(const.)}{n^{c_-}}, \tag{2.13}$$

for some positive constant $c_-$. Notice that by Proposition 2.2 we have $\mathcal{Z}^{\mathrm{p}}_{0,N} = \varphi^{(0,0)}_{N+1}(0,0)$ and $\mathcal{Z}^{\mathrm{w}}_{0,N} = \mathcal{Z}^{\mathrm{p}}_{0,N} \cdot w_{0,0}(N+1)$, hence by (2.10) and (2.13) we have for every $N \in \mathbb{N}$

$$\mathcal{Z}^{\mathrm{p}}_{0,N} \geq \mathcal{Z}^{\mathrm{w}}_{0,N} \geq \frac{(const.)}{N^{2+c_-}}, \tag{2.14}$$

so that the last inequality in (1.9) is proven.

We prove the lower bound (2.13) in Appendix C.1; the idea is to restrict the expectation that defines $w_{0,0}(n)$ on a suitable subset of paths, whose probability can be estimated. On the other hand, the upper bound (2.12) follows directly combining the following lemma with the upper bound in (1.23) (which is proven in Appendix C.2).

LEMMA 2.4. *There exists a positive constant $C$ such that for every $n \in \mathbb{N}$*

$$\sup_{0\leq x,y\leq \sqrt{n}} w_{x,y}(n) \leq C \cdot \mathbf{P}^{(0,0)}(\Omega^+_{\lfloor n/2 \rfloor}). \tag{2.15}$$

PROOF. The inclusion bound gives (for $x, y \geq 0$)

$$\begin{aligned} w_{x,y}(n) &= \mathbf{P}^{(-x,0)}(\Omega^+_{n-2} | Z_{n-1} = y, Z_n = 0) \\ &\leq \mathbf{P}^{(-x,0)}(\Omega^+_{\lfloor n/2 \rfloor} | Z_{n-1} = y, Z_n = 0), \end{aligned}$$

and disintegrating over the variables $Z_{\lfloor n/2\rfloor -1}, Z_{\lfloor n/2 \rfloor}$ we get [recall (2.6) and (2.9)]

$$\begin{aligned} &w_{x,y}(n) \\ &\quad \leq \int_{a,b\in\mathbb{R}^+} \frac{\mathbf{P}^{(-x,0)}(\Omega^+_{\lfloor n/2\rfloor}, Z_{\lfloor n/2\rfloor -1} \in da, Z_{\lfloor n/2 \rfloor} \in db) \cdot \varphi^{(b-a,b)}_{n-\lfloor n/2\rfloor}(-y,0)}{\varphi^{(-x,0)}_n(-y,0)}. \end{aligned}$$

Proposition 2.3 together with (2.9) yields

$$\sup_{a,b,y\in\mathbb{R}} \varphi^{(b-a,b)}_{n-\lfloor n/2\rfloor}(-y,0) \le \frac{(const.)}{n^2}, \qquad \inf_{0\le x,y\le\sqrt{n}} \varphi^{(-x,0)}_{n}(-y,0) \ge \frac{(const.')}{n^2},$$

hence setting $C := (const.)/(const.') > 0$ we have

$$\begin{aligned} w_{x,y}(n) &\le C \int_{a,b\in\mathbb{R}^+} \mathbf{P}^{(-x,0)}(\Omega^+_{\lfloor n/2\rfloor}, Z_{\lfloor n/2\rfloor -1} \in da, Z_{\lfloor n/2\rfloor} \in db) \\ &= C\mathbf{P}^{(-x,0)}(\Omega^+_{\lfloor n/2\rfloor}), \end{aligned}$$

and the proof is completed noting that $\mathbf{P}^{(-x,0)}(\Omega^+_{\lfloor n/2\rfloor}) \le \mathbf{P}^{(0,0)}(\Omega^+_{\lfloor n/2\rfloor})$ by (2.3). ☐

**3. The interacting case $\varepsilon > 0$: a renewal theory description.** In this section we study in detail the laws $\mathbb{P}^{\mathrm{p}}_{\varepsilon,N}$ and $\mathbb{P}^{\mathrm{w}}_{\varepsilon,N}$ in the case $\varepsilon > 0$. The crucial result is that the contact set $\{i \in \mathbb{Z}^+ : \varphi_i = 0\}$ can be described in terms of a Markov renewal process. Throughout the section we assume that $\varepsilon > 0$.

3.1. *The law of the contact set.* We introduce the *contact set* $\tau$ by:

$$\tau := \{i \in \mathbb{Z}^+ : \varphi_i = 0\} \subset \mathbb{Z}^+, \tag{3.1}$$

where we set by definition $\varphi_0 = 0$, so that $0 \in \tau$. It is practical to identify the set $\tau$ with the increasing sequence of random variables $\{\tau_k\}_{k\ge 0}$ defined by

$$\tau_0 := 0, \qquad \tau_{k+1} := \inf\{i > \tau_k : \varphi_i = 0\}. \tag{3.2}$$

Observe that the random variable $\ell_N$, introduced in (1.12), may be expressed as $\ell_N := \max\{k : \tau_k \le N\}$. Next we introduce the process $\{J_k\}_{k\ge0}$ that gives the height of the field before the contact points:

$$J_0 := 0, \qquad J_k := \varphi_{\tau_k - 1}. \tag{3.3}$$

The basic observation is that the joint law of the process $\{\ell_N, (\tau_k)_{k\le \ell_N}, (J_k)_{k\le\ell_N}\}$ under $\mathbb{P}^a_{\varepsilon,N}$ can be written in the following "product form": for $k \in \mathbb{N}$, for $(t_i)_{i=1,\ldots,k} \in \mathbb{N}$ with $0 < t_1 < \cdots < t_{k-1} < t_k := N$ and for $(y_i)_{i=1,\ldots,k} \in \mathbb{R}$ we have

$$\begin{aligned} &\mathbb{P}^a_{\varepsilon,N}(\ell_N = k, \tau_i = t_i, J_i \in dy_i, i = 1,\ldots,k) \\ &\quad = \frac{1}{\mathcal{Z}^a_{\varepsilon,N}} \, \varepsilon^{k-1} \mathsf{F}^a_{0,dy_1}(t_1) \\ &\qquad \times \mathsf{F}^a_{y_1,dy_2}(t_2 - t_1) \cdot \cdots \cdot \mathsf{F}^a_{y_{k-1},dy_k}(N - t_{k-1}) \cdot \mathsf{F}^a_{y_k,\{0\}}(1), \end{aligned} \tag{3.4}$$

for a suitable kernel $\mathsf{F}^a_{x,dy}(n)$ that we now define. For $a = \mathrm{p}$ we have

$$\mathsf{F}^{\mathrm{p}}_{x,dy}(n) := \begin{cases} e^{-\mathcal{H}_{[-1,1]}(x,0,0)}\delta_0(dy) = e^{-V(x)}\delta_0(dy), & \text{if } n=1, \\ e^{-\mathcal{H}_{[-1,2]}(x,0,y,0)}\,dy = e^{-V(x+y)-V(-2y)}\,dy, & \text{if } n=2, \\ \left\{ \begin{array}{l} \left(\int_{\mathbb{R}^{n-2}} e^{-\mathcal{H}_{[-1,n]}(\varphi_{-1},\dots,\varphi_n)}\,d\varphi_1\cdots d\varphi_{n-2}\right) dy \\ \text{where } \varphi_{-1}=x, \varphi_0=0, \varphi_{n-1}=y, \varphi_n=0 \end{array} \right\}, & \text{if } n \geq 3, \end{cases} \tag{3.5}$$

and the definition of $\mathsf{F}^{\mathrm{w}}_{x,dy}(n)$ is analogous: we just have to impose that $x, y \geq 0$ and for $n \geq 3$ we also have to restrict the integral in (3.5) on $(\mathbb{R}^+)^{n-2}$. Although these formulas may appear quite involved, they follow easily from the definition of $\mathbb{P}^a_{\varepsilon,N}$. In fact it suffice, to expand the product of measures in the r.h.s. of (1.5) and (1.7) as a sum of "monomials," according to the elementary formula (where we set $\Omega^{\mathrm{p}} := \mathbb{R}$ and $\Omega^{\mathrm{w}} := \mathbb{R}^+$)

$$\prod_{i=1}^{N-1}(\varepsilon\delta_0(d\varphi_i) + d\varphi_i \mathbf{1}_{(\varphi_i\in\Omega^a)}) = \sum_{A\subset\{1,\dots,N-1\}} \varepsilon^{|A|} \prod_{m\in A}\delta_0(d\varphi_m) \prod_{n\in A^{\complement}} d\varphi_n \mathbf{1}_{(\varphi_n\in\Omega^a)}.$$

It is then clear that $A = \{\tau_1, \dots, \tau_{\ell_N - 1}\}$ and integrating over the variables $\varphi_i$ with index $i \notin A \cup (A-1)$ one gets to (3.4). We stress that the algebraic structure of (3.4) retains a crucial importance, that we are going to exploit in the next paragraph.

From (3.5) it follows that the kernel $\mathsf{F}^a_{x,dy}(n)$ is a Dirac mass in $y$ for $n = 1$ while it is absolutely continuous for $n \geq 2$. Then it is convenient to introduce the $\sigma$-finite measure $\mu(dx) := \delta_0(dx) + dx$, so that we can write

$$\mathsf{F}^a_{x,dy}(n) = \mathsf{f}^a_{x,y}(n)\,\mu(dy). \tag{3.6}$$

The interesting fact is that the density $\mathsf{f}^a_{x,y}(n)$ can be rephrased explicitly in terms of the process $\{Z_k\}_k$ introduced in Section 2.1. Let us start with the pinning case $a = \mathrm{p}$: From Lemma 2.1 and from (3.5) it follows that, for $n \geq 2$, $\mathsf{f}^{\mathrm{p}}_{x,y}(n)$ is nothing but the density of $(Z_{n-1}, Z_n)$ at $(y, 0)$, under $\mathbf{P}^{(-x,0)}$. Recalling the definition (2.8) of $\varphi^{(\cdot,\cdot)}_n(\cdot,\cdot)$ and the fact that $Z_{n-1} = Z_n - Y_n$, we can write

$$\mathsf{f}^{\mathrm{p}}_{x,y}(n) = \begin{cases} e^{-V(x)}\mathbf{1}_{(y=0)}, & \text{if } n=1, \\ \varphi^{(-x,0)}_n(-y,0)\mathbf{1}_{(y\neq 0)} & \\ \qquad = \varphi^{(0,0)}_n(-y+x, nx)\mathbf{1}_{(y\neq0)}, & \text{if } n\geq 2, \end{cases} \tag{3.7}$$

where we have used relation (2.9). Analogously, recalling the definition (2.11) of $w_{x,y}(n)$, in the wetting case we have

$$\mathsf{f}^{\mathrm{w}}_{x,y}(n) = \mathsf{f}^{\mathrm{p}}_{x,y}(n) \cdot w_{x,y}(n). \tag{3.8}$$

Equations (3.6), (3.7) and (3.8) provide a description of the kernel $\mathsf{F}^a_{x,dy}(n)$ which is both simpler and more useful than the original definition (3.5).

3.2. *A Markov renewal theory interpretation.* Equation (3.4) expresses the law of $\{(\tau_k, J_k)\}_k$ under $\mathbb{P}^a_{\varepsilon,N}$ in terms of an explicit kernel $\mathsf{F}^a_{x,dy}(n)$. The crucial point is that the algebraic structure of equation (3.4) allows to modify the kernel, in order to give this formula a direct renewal theory interpretation. In fact we set

$$\mathsf{K}^{a,\varepsilon}_{x,dy}(n) := \varepsilon \mathsf{F}^a_{x,dy}(n) e^{-\mathrm{F}^a(\varepsilon) n} \frac{v^a_\varepsilon(y)}{v^a_\varepsilon(x)}, \tag{3.9}$$

where the number $\mathrm{F}^a(\varepsilon) \in [0,\infty)$ and the positive real function $v^a_\varepsilon(\cdot)$ will be defined explicitly in Section 4. Of course this is an abuse of notation, because the symbol $\mathrm{F}^a(\varepsilon)$ was already introduced to denote the free energy [cf. (1.8)], but we will show in Section 5.2 that the two quantities indeed coincide. We denote by $\mathsf{k}^{a,\varepsilon}_{x,y}(n)$ the density of $\mathsf{K}^{a,\varepsilon}_{x,dy}(n)$ with respect to $\mu(dy)$, that is,

$$\mathsf{k}^{a,\varepsilon}_{x,y}(n) := \varepsilon \mathsf{f}^a_{x,y}(n) e^{-\mathrm{F}^a(\varepsilon) n} \frac{v^a_\varepsilon(y)}{v^a_\varepsilon(x)}. \tag{3.10}$$

The reason for introducing the kernel $\mathsf{K}^{a,\varepsilon}_{x,dy}(n)$ lies in the following fundamental fact: the number $\mathrm{F}^a(\varepsilon)$ and the function $v^a_\varepsilon(\cdot)$ appearing in (3.9) can be chosen such that,

$$\forall x \in \mathbb{R}: \qquad \int_{y\in\mathbb{R}} \sum_{n\in\mathbb{N}} \mathsf{K}^{a,\varepsilon}_{x,dy}(n) = \frac{\varepsilon}{\varepsilon^a_c} \wedge 1 \le 1\,, \tag{3.11}$$

where $\varepsilon^a_c \in (0,\infty)$ is a fixed number. A detailed proof and discussion of this fact, with an explicit definition of $\varepsilon^a_c$, $\mathrm{F}^a(\varepsilon)$ and $v^a_\varepsilon(\cdot)$, is deferred to Section 4; for the moment we focus on its consequences.

Thanks to (3.11), we can define the law $\mathcal{P}^a_\varepsilon$ under which the joint process $\{(\tau_k, J_k)\}_{k\ge 0}$ is a (possibly defective) Markov chain on $\mathbb{Z}^+ \times \mathbb{R}$, with starting value $(\tau_0, J_0) = (0,0)$ and with transition kernel given by

$$\mathcal{P}^a_\varepsilon((\tau_{k+1}, J_{k+1}) \in (\{n\}, dy) | (\tau_k, J_k) = (m,x)) = \mathsf{K}^{a,\varepsilon}_{x,dy}(n-m). \tag{3.12}$$

An alternative (and perhaps more intuitive) definition is as follows:

- First sample the process $\{J_k\}_{k\ge 0}$ as a (defective if $\varepsilon < \varepsilon^a_c$) Markov chain on $\mathbb{R}$, with $J_0 = 0$ and with transition kernel

$$\mathcal{P}^a_\varepsilon(J_{k+1} \in dy \mid J_k = x) = \sum_{n\in\mathbb{N}} \mathsf{K}^{a,\varepsilon}_{x,dy}(n) =: D^{a,\varepsilon}_{x,dy}. \tag{3.13}$$

  In the defective case we take $\infty$ as cemetery.
- Then sample the increments $\{T_k := \tau_k - \tau_{k-1}\}_{k\in\mathbb{N}}$ as a sequence of independent, but not identically distributed, random variables, according to the conditional law:

$$\mathcal{P}^a_\varepsilon(T_k = n \,|\, \{J_i\}_{i\ge 0}) = \begin{cases} \mathbf{1}_{(n=1)}, & \text{if } J_k = 0, \\ \dfrac{\mathsf{k}^{a,\varepsilon}_{J_{k-1},J_k}(n) \mathbf{1}_{(n\ge 2)}}{\sum_{m\ge 2} \mathsf{k}^{a,\varepsilon}_{J_{k-1},J_k}(m)}, & \text{if } J_k \ne 0,\ J_k \ne \infty, \\ \mathbf{1}_{(n=\infty)} & \text{if } J_k = \infty. \end{cases}$$

We stress that the process $\{(\tau_k, J_k)\}_{k\geq 0}$ is defective if $\varepsilon < \varepsilon_c^a$ and proper if $\varepsilon \geq \varepsilon_c^a$; cf. (3.11). The process $\{\tau_k\}_{k\geq 0}$ under $\mathcal{P}_\varepsilon^a$ is what is called a *Markov renewal process* and $\{J_k\}_{k\geq 0}$ is its *modulating chain*. This is a generalization of classical renewal processes, because $\tau_n = T_1 + \cdots + T_n$ where the variables $\{T_k\}_{k\in\mathbb{N}}$ are allowed to have a special kind of dependence, namely they are independent conditionally on the modulating chain $\{J_k\}_{k\geq 0}$. For a detailed account on Markov renewal processes we refer to [4].

Now let us come back to (3.4). We perform the substitution $\mathsf{F}^a_{x,dy}(n) \to \mathsf{K}^{a,\varepsilon}_{x,dy}(n)$, defined in (3.9): the boundary terms $v_\varepsilon^a(y)/v_\varepsilon^a(x)$ get simplified and the exponential term $e^{-\mathrm{F}^a(\varepsilon)n}$ factorizes, so that we get

$$\begin{aligned}
&\mathbb{P}^a_{\varepsilon,N}(\ell_N = k, \tau_i = t_i, J_i \in dy_i, i = 1, \dots, k) \\
&\qquad = \frac{e^{\mathrm{F}^a(\varepsilon)N}}{\varepsilon^2\, \mathcal{Z}^a_{\varepsilon,N}} \mathsf{K}^{a,\varepsilon}_{0,dy_1}(t_1) \cdot \mathsf{K}^{a,\varepsilon}_{y_1,dy_2}(t_2 - t_1) \cdots \\
&\qquad\qquad \cdot \mathsf{K}^{a,\varepsilon}_{y_{k-1},dy_k}(N - t_{k-1}) \cdot \mathsf{K}^{a,\varepsilon}_{y_k,\{0\}}(1).
\end{aligned} \tag{3.14}$$

Moreover, since the partition function $\mathcal{Z}^a_{\varepsilon,N}$ is the normalizing constant that makes $\mathbb{P}^a_{\varepsilon,N}$ a probability, it can be expressed as

$$\begin{aligned}
\mathcal{Z}^a_{\varepsilon,N} = \frac{e^{\mathrm{F}^a(\varepsilon)N}}{\varepsilon^2} \sum_{k=1}^{N} \sum_{\substack{t_i\in\mathbb{N},\, i=1,\dots,k \\ 0<t_1<\cdots<t_k:=N}} \int_{\mathbb{R}^k} \left( \prod_{i=1}^{k} \mathsf{K}^{a,\varepsilon}_{y_{i-1},dy_i}(t_i - t_{i-1}) \right) \\
\cdot \mathsf{K}^{a,\varepsilon}_{y_k,\{0\}}(1).
\end{aligned} \tag{3.15}$$

We are finally ready to make explicit the link between the law $\mathcal{P}_\varepsilon^a$ and our model $\mathbb{P}^a_{\varepsilon,N}$. Let us introduce the event

$$\mathcal{A}_N := \{\, \{N, N+1\} \subset \tau \} = \{\exists k \geq 0 : \tau_k = N,\, \tau_{k+1} = N+1\}. \tag{3.16}$$

The following proposition is an immediate consequence of (3.12), (3.14) and (3.15).

Proposition 3.1. *For any $N \in \mathbb{N}$ and $\varepsilon > 0$, the vector $\{\ell_N, (\tau_i)_{i\leq \ell_N}, (J_i)_{i\leq \ell_N}\}$ has the same law under $\mathbb{P}^a_{\varepsilon,N}$ and under the conditional law $\mathcal{P}_\varepsilon^a(\cdot|\mathcal{A}_N)$: for all $k, \{t_i\}_i$ and $\{y_i\}_i$*

$$\begin{aligned}
&\mathbb{P}^a_{\varepsilon,N}(\ell_N = k, \tau_i = t_i, J_i \in dy_i, i \leq k) \\
&\qquad = \mathcal{P}^a_\varepsilon(\ell_N = k, \tau_i = t_i, J_i \in dy_i, i \leq k \mid \mathcal{A}_N).
\end{aligned}$$

*Moreover the partition function can be expressed as $\mathcal{Z}^a_{\varepsilon,N} = (e^{\mathrm{F}^a(\varepsilon)N}/\varepsilon^2) \cdot \mathcal{P}^a_\varepsilon(\mathcal{A}_N)$.*

Thus we have shown that the contact set $\tau \cap [0, N]$ under the pinning law $\mathbb{P}^a_{\varepsilon,N}$ is distributed like a Markov renewal process (of law $\mathcal{P}^a_\varepsilon$) conditioned to visit $\{N, N+1\}$. The crucial point is that $\mathcal{P}^a_\varepsilon$ does not have any dependence on $N$, therefore all the dependence on $N$ of $\mathbb{P}^a_{\varepsilon,N}$ is contained in the conditioning on the event $\mathcal{A}_N$. As will be clear in the next sections, this fact is the key to all our results.

**4. An infinite-dimensional Perron–Frobenius problem.** In this section we prove that for every $\varepsilon > 0$ the nonnegative number $\mathrm{F}^a(\varepsilon)$ and the positive real function $v^a_\varepsilon(\cdot) : \mathbb{R} \to (0, \infty)$ appearing in the definition of $\mathsf{K}^{a,\varepsilon}_{x,dy}(n)$, [cf. (3.9)], can be chosen in such a way that (3.11) holds true.

4.1. *Some analytical preliminaries.* We recall that the kernel $\mathsf{F}^a_{x,dy}(n)$ and its density $\mathsf{f}^a_{x,y}(n)$ are defined in (3.6), (3.7) and (3.8), and that $\mu(dx) = \delta_0(dx) + dx$. In particular we have $0 \le \mathsf{f}^{\mathrm{w}}_{x,y}(n) \le \mathsf{f}^{\mathrm{p}}_{x,y}(n)$. We first list some important properties of $\mathsf{f}^{\mathrm{p}}_{x,y}(n)$:

- Uniformly for $x, y$ in compact sets we have

$$\mathsf{f}^{\mathrm{p}}_{x,y}(n) \sim \frac{c}{n^2} \qquad (n \to \infty), \tag{4.1}$$

where $c := 6/(\pi\sigma^2)$; moreover there exists $C > 0$ such that $\forall x, y \in \mathbb{R}$ and $n \in \mathbb{N}$

$$\mathsf{f}^{\mathrm{p}}_{x,y}(n) \le \frac{C}{n^2}, \qquad \int_{\mathbb{R}} dz\, \mathsf{f}^{\mathrm{p}}_{x,z}(n) \le \frac{C}{n^{3/2}}, \qquad \int_{\mathbb{R}} dz\, \mathsf{f}^{\mathrm{p}}_{z,y}(n) \le \frac{C}{n^{3/2}}. \tag{4.2}$$

Both the above relations follow comparing (3.7) with Proposition (2.3).
- For $n \ge 2$ we have

$$\int_{\mathbb{R}^2} dx\, dy\, \mathsf{f}^{\mathrm{p}}_{x,y}(n) = \frac{1}{n}, \tag{4.3}$$

as follows from (3.7) recalling that $\varphi^{(0,0)}_n(\cdot,\cdot)$ is a probability density.

For $\lambda \ge 0$ we introduce the kernel

$$B^{a,\lambda}_{x,dy} := \sum_{n\in\mathbb{N}} e^{-\lambda n} \mathsf{F}^a_{x,dy}(n), \tag{4.4}$$

which induces the integral operator: $(B^{a,\lambda} h)(x) := \int_{\mathbb{R}} B^{a,\lambda}_{x,dy}\, h(y)$. Note that for every $x \in \mathbb{R}$ the kernel $B^{a,\lambda}_{x,dy}$ is absolutely continuous with respect to the measure $\mu$, so that we can write $B^{a,\lambda}_{x,dy} = b^{a,\lambda}(x,y)\, \mu(dy)$, where the density $b^{a,\lambda}(x,y)$ is given by [cf. (3.7)]

$$b^{a,\lambda}(x,y) = e^{-\lambda}\, \mathsf{f}^a_{x,0}(1) \mathbf{1}_{(y=0)} + \sum_{n\ge 2} e^{-\lambda n}\, \mathsf{f}^a_{x,y}(n) \mathbf{1}_{(y\neq 0)}. \tag{4.5}$$

The following result is of basic importance.

LEMMA 4.1. *For every $\lambda \geq 0$, $B^{a,\lambda}$ is a compact operator on the Hilbert space $L^2(\mathbb{R}, d\mu)$.*

PROOF. We are going to check the stronger condition that $B^{a,\lambda}$ is Hilbert–Schmidt, that is,

$$\int_{\mathbb{R}^2} b^{a,\lambda}(x,y)^2 \mu(dx)\mu(dy) < \infty. \tag{4.6}$$

Since $0 \leq \mathsf{f}^{\mathrm{w}}_{x,y}(n) \leq \mathsf{f}^{\mathrm{p}}_{x,y}(n)$ it suffices to focus on the case $a = \mathrm{p}$. Setting $\lambda = 0$ in the r.h.s. of (4.5) we obtain

$$b^{\mathrm{p},\lambda}(x,y)^2 \leq \mathsf{f}^{\mathrm{p}}_{x,0}(1)^2 \mathbf{1}_{(y=0)} + \sum_{n,m\geq 2} \mathsf{f}^{\mathrm{p}}_{x,y}(n)\mathsf{f}^{\mathrm{p}}_{x,y}(m)\mathbf{1}_{(y\neq 0)},$$

hence

$$\begin{aligned} &\int_{\mathbb{R}^2} b^{\mathrm{p},\lambda}(x,y)^2 \mu(dx)\mu(dy) \\ &\qquad \leq \int_{\mathbb{R}} \mathsf{f}^{\mathrm{p}}_{x,0}(1)^2 \mu(dx) + \sum_{n,m\geq 2} \int_{\mathbb{R}} \mathsf{f}^{\mathrm{p}}_{0,y}(n)\mathsf{f}^{\mathrm{p}}_{0,y}(m)\, dy \\ &\qquad\quad + \sum_{n,m\geq 2} \int_{\mathbb{R}^2} \mathsf{f}^{\mathrm{p}}_{x,y}(n)\mathsf{f}^{\mathrm{p}}_{x,y}(m)\, dx\, dy. \end{aligned} \tag{4.7}$$

From the definition (3.7) of $\mathsf{f}^{\mathrm{p}}_{x,0}(1)$ it is immediate to check that the first term in the r.h.s. is finite. For the second term, we apply the first and the second relation in (4.2), getting

$$\sum_{n,m\geq 2} \int_{\mathbb{R}} \mathsf{f}^{\mathrm{p}}_{0,y}(n)\mathsf{f}^{\mathrm{p}}_{0,y}(m)\, dy \leq \sum_{n,m\geq 2} \frac{C}{n^2} \int_{\mathbb{R}} \mathsf{f}^{\mathrm{p}}_{0,y}(m)\, dy \leq \sum_{n,m\geq 2} \frac{C}{n^2} \frac{C}{m^{3/2}} < \infty.$$

For the third term, it is convenient first to exploit the symmetry between $n$ and $m$, restricting the sum on the set $n \geq m$:

$$\sum_{n,m\geq 2} \int_{\mathbb{R}^2} \mathsf{f}^{\mathrm{p}}_{x,y}(n)\mathsf{f}^{\mathrm{p}}_{x,y}(m)\, dx\, dy \leq 2 \sum_{m\geq 2} \sum_{n\geq m} \int_{\mathbb{R}^2} \mathsf{f}^{\mathrm{p}}_{x,y}(n)\mathsf{f}^{\mathrm{p}}_{x,y}(m)\, dx\, dy.$$

Then applying relations (4.2) and (4.3) we get

$$\sum_{m\geq 2} \sum_{n\geq m} \int_{\mathbb{R}^2} \mathsf{f}^{\mathrm{p}}_{x,y}(n)\mathsf{f}^{\mathrm{p}}_{x,y}(m)\, dx\, dy \leq \sum_{m\geq 2} \sum_{n\geq m} \frac{1}{m} \frac{C}{n^2} = C \sum_{m\geq 2} \frac{1}{m} \sum_{n\geq m} \frac{1}{n^2}.$$

However, the last sum is bounded by $(const.)/m$, hence the r.h.s. is finite. □

4.2. *A formula for the free energy.* Lemma 4.1 allows us to apply an infinite-dimensional analogue of the classical Perron–Frobenius theorem. We first introduce the function $\delta^a(\lambda) \in [0,\infty)$ defined for $\lambda \geq 0$ by

$$\delta^a(\lambda) := \text{spectral radius of the operator } B^{a,\lambda}. \tag{4.8}$$

Notice that $\{0\}$ is a *proper atom*, hence a *small set*, of the kernel $B^{a,\lambda}$; cf. [26], Section 4.2. Therefore by [26], Section 3.2 one can define $\delta^a(\lambda)$ equivalently as

$$\delta^a(\lambda) := \inf\left\{ R > 0 : \sum_{n=0}^{\infty} \frac{(B^{a,\lambda})^{\circ n}_{0,\{0\}}}{R^n} < \infty \right\}, \tag{4.9}$$

where the convolution $\circ$ between kernels is defined in Section 1.7. One checks easily that indeed $\delta^a(\lambda) \in (0,\infty)$ for every $\lambda \geq 0$. By Theorem 1 in [33], $\delta^a(\lambda)$ is an isolated and *simple* eigenvalue of $B^{a,\lambda}$. The function $\delta^a(\cdot)$ is nonincreasing, continuous on $[0,\infty)$ and analytic on $(0,\infty)$, because the operator $B^{a,\lambda}$ has these properties and $\delta^a(\lambda)$ is simple and isolated; cf. [20], Chapter VII, Section 1.3. The analyticity and the fact that $\delta^a(\cdot)$ is not constant (because $\delta^a(\lambda) \to 0$ as $\lambda \to +\infty$) force $\delta^a(\cdot)$ to be strictly decreasing. Denoting by $(\delta^a)^{-1}(\cdot)$ the inverse function, defined on the domain $(0, \delta^a(0)]$, we can now define $\varepsilon_c^a$ and $\mathrm{F}^a(\varepsilon)$ by

$$\varepsilon_c^a := \frac{1}{\delta^a(0)}, \qquad \mathrm{F}^a(\varepsilon) := \begin{cases} (\delta^a)^{-1}(1/\varepsilon), & \text{if } \varepsilon \geq \varepsilon_c^a, \\ 0, & \text{if } \varepsilon \leq \varepsilon_c^a. \end{cases} \tag{4.10}$$

From now on we focus on the operator $B^{a,\mathrm{F}^a(\varepsilon)}$, that is, on $B^{a,\lambda}$ for $\lambda = \mathrm{F}^a(\varepsilon)$, whose spectral radius equals $1/\varepsilon \wedge 1/\varepsilon_c^a$ by construction. Notice that $b^{a,\mathrm{F}^a(\varepsilon)}(x,y) > 0$ for every $x, y \in \Omega^a$, where $\Omega^{\mathrm{p}} = \mathbb{R}$ and $\Omega^{\mathrm{w}} = \mathbb{R}^+$. Then Theorem 1 in [33] gives the existence of the so-called right and left *Perron–Frobenius eigenfunctions* of $B^{a,\mathrm{F}^a(\varepsilon)}$, that is, of two functions $v_\varepsilon^a(\cdot), w_\varepsilon^a(\cdot) \in L^2(\mathbb{R}, d\mu)$ such that $v_\varepsilon^a(x) > 0$ and $w_\varepsilon^a(x) > 0$, for $\mu$-a.e. $x \in \Omega^a$, and such that

$$\begin{aligned} \int_{y\in\mathbb{R}} B^{a,\mathrm{F}^a(\varepsilon)}_{x,dy} v_\varepsilon^a(y) &= \left(\frac{1}{\varepsilon} \wedge \frac{1}{\varepsilon_c^a}\right) v_\varepsilon^a(x), \\ \int_{x\in\mathbb{R}} w_\varepsilon^a(x) B^{a,\mathrm{F}^a(\varepsilon)}_{x,dy} \mu(dx) &= \left(\frac{1}{\varepsilon} \wedge \frac{1}{\varepsilon_c^a}\right) w_\varepsilon^a(y)\mu(dy). \end{aligned} \tag{4.11}$$

Notice that in the wetting case $v_\varepsilon^{\mathrm{w}}(x) = w_\varepsilon^{\mathrm{w}}(x) = 0$ for all $x < 0$. Spelling out the first equation in (4.11) we have

$$v_\varepsilon^a(x) = \frac{1}{1/\varepsilon \wedge 1/\varepsilon_c^a} \sum_{n\in\mathbb{N}} e^{-\mathrm{F}^a(\varepsilon) n} \int_{y\in\mathbb{R}} \mathsf{f}^a_{x,y}(n)\, v_\varepsilon^a(y)\, \mu(dy). \tag{4.12}$$

This yields easily that $v_\varepsilon^a(x) > 0$ for every $x \in \Omega^a$ (and not only $\mu$-a.e.). One shows analogously that $w_\varepsilon^a(x) > 0$ for every $x \in \Omega^a$.

Having defined the quantities $\varepsilon_c^a, \mathrm{F}^a(\varepsilon)$ and $v_\varepsilon^a(\cdot)$, it remains to check that (3.11) indeed holds true. But this is a straightforward consequence of (4.11); in fact by the definition (3.9) of $\mathsf{K}_{x,dy}^{a,\varepsilon}(n)$ we have

$$\int_{y\in\mathbb{R}} \sum_{n\in\mathbb{N}} \mathsf{K}_{x,dy}^{a,\varepsilon}(n) = \frac{\varepsilon}{v_\varepsilon^a(x)} \int_{y\in\mathbb{R}} \left( \sum_{n\in\mathbb{N}} \mathsf{F}_{x,dy}^a(n) e^{-\mathrm{F}^a(\varepsilon)n} \right) v_\varepsilon^a(y)$$
$$= \frac{\varepsilon}{v_\varepsilon^a(x)} \int_{y\in\mathbb{R}} B_{x,dy}^{a,\mathrm{F}^a(\varepsilon)} v_\varepsilon^a(y) = \varepsilon \left( \frac{1}{\varepsilon} \wedge \frac{1}{\varepsilon_c^a} \right) = 1 \wedge \frac{\varepsilon}{\varepsilon_c^a}.$$

We note that the functions $v_\varepsilon^a(\cdot)$ and $w_\varepsilon^a(\cdot)$ are uniquely defined up to a multiplicative constant and we use this degree of freedom to fix $\langle v_\varepsilon^a, w_\varepsilon^a \rangle = \int_{\mathbb{R}} v_\varepsilon^a w_\varepsilon^a \, d\mu = 1$. In other words, the measure $\nu_\varepsilon^a$ defined by

$$\nu_\varepsilon^a(dx) := w_\varepsilon^a(x) v_\varepsilon^a(x) \mu(dx) \tag{4.13}$$

is a probability measure: $\nu_\varepsilon^a(\mathbb{R}) = 1$. An important observation is that for $\varepsilon \geq \varepsilon_c^a$ the probability measure $\nu_\varepsilon^a$ is invariant for the kernel $D_{x,dy}^{a,\varepsilon}$, as follows from (3.13) and (3.9). This means that the Markov chain $\{J_k\}_k$ is *positive recurrent*; cf. [26].

We conclude this section with a simple perturbation result that will be useful later. Let $A_{x,\,dy}$ and $C_{x,\,dy}$ be two nonnegative kernels that induce two compact operators on $L^2(\mathbb{R}, d\mu)$. Assume that the spectral radius of $C_{x,\,dy}$ is strictly positive. For $t \in [0, \infty)$ we set $\gamma(t) :=$ spectral radius of $(A_{x,\,dy} + t \cdot C_{x,\,dy})$. Then we have the following:

LEMMA 4.2. *The function $t \mapsto \gamma(t)$ is strictly increasing; in particular $\gamma(0) < \gamma(1)$.*

PROOF. The function $\gamma(\cdot)$ is clearly nondecreasing. It follows by Theorem 1 in [33] that $\gamma(t)$ is a simple and isolated eigenvalue of $A_{x,dy} + t \cdot C_{x,dy}$ for every $t \geq 0$; therefore by perturbation theory [20], Chapter VII, Section 1.3 the function $\gamma(\cdot)$ is analytic. Since $\gamma(t) \geq (const.) \cdot t$ (here we use the hypothesis that the spectral radius of $C_{x,dy}$ is strictly positive), the function $\gamma(\cdot)$ is nonconstant and therefore it must be strictly increasing. □

**5. Proof of Theorem 1.2.** We apply the content of the preceding sections to prove Theorem 1.2. Before that, we show that from the Markov renewal structure described in Section 3 one can extract a genuine renewal process, which will be a basic technical tool.

5.1. *From Markov renewals to true renewals.* Recalling the definition (3.1) of the contact set $\tau$, we introduce the subset $\chi$ of the *adjacent contact points* defined by

$$\chi := \{ i \in \mathbb{Z}^+ : \varphi_{i-1} = \varphi_i = 0 \} \subset \tau \subset \mathbb{Z}^+ , \tag{5.1}$$

where we set by definition $\varphi_{-1} = \varphi_0 = 0$, so that $0 \in \chi$. Note that $\chi = \tau \cap (\tau + 1)$ and that the points of $\chi$ are the random variables $\{\chi_k\}_{k \geq 0}$ defined by

$$\chi_0 := 0, \qquad \chi_{k+1} := \inf\{i > \chi_k : \varphi_{i-1} = \varphi_i = 0\}. \tag{5.2}$$

By (3.16) the event $\mathcal{A}_N$ can be written as $\mathcal{A}_N = \{N + 1 \in \chi\}$; therefore by Proposition 3.1 the partition function can be written as

$$\mathcal{Z}^a_{\varepsilon,N} = \frac{e^{\mathrm{F}^a(\varepsilon)N}}{\varepsilon^2} \cdot \mathcal{P}^a_\varepsilon(N + 1 \in \chi). \tag{5.3}$$

The reason for focusing on the process $\{\chi_k\}_k$ is explained by the following proposition.

Proposition 5.1. *The process $\{\chi_k\}_{k \geq 0}$ under $\mathcal{P}^a_\varepsilon$ is a classical (i.e., not Markov) renewal process, which is nonterminating for $\varepsilon \geq \varepsilon^a_c$.*

Proof. We introduce the epochs $\{\zeta_k\}_k$ of return to zero of the process $\{J_k\}_k$:

$$\zeta_0 := 0, \qquad \zeta_{n+1} := \inf\{k > \zeta_n : J_k = 0\}. \tag{5.4}$$

It is clear that the variables $\{\zeta_k - \zeta_{k-1}\}_{k \in \mathbb{N}}$ are i.i.d. under $\mathcal{P}^a_\varepsilon$, because $\{J_k\}_k$ is a Markov chain. Observe that $\chi_k = \tau_{\zeta_k}$, for every $k \geq 0$, as it follows by (3.2), (3.3) and (5.2). This fact implies that also the variables $\{\chi_k - \chi_{k-1}\}_{k \in \mathbb{N}} = \{\tau_{\zeta_k} - \tau_{\zeta_{k-1}}\}_{k \in \mathbb{N}}$ are i.i.d. under $\mathcal{P}^a_\varepsilon$, because the transition kernel of the process $\{(\tau_k, J_k)\}_{k \geq 0}$ is a function of $(n - m)$; cf. the r.h.s. of (3.12). Therefore $\{\chi_k\}_{k \geq 0}$ under $\mathcal{P}^a_\varepsilon$ is a genuine renewal process.

We have already observed that for $\varepsilon \geq \varepsilon^a_c$ the Markov chain $\{J_k\}_{k \geq 0}$ is *positive recurrent*, because for $\varepsilon \geq \varepsilon^a_c$ the probability measure $\nu^a_\varepsilon$ defined in (4.13) is by construction an invariant measure for its transition kernel $D^{a,\varepsilon}_{x,dy}$; cf. (3.13) and (3.9). Since $\nu^a_\varepsilon(\{0\}) > 0$, the state 0 is an *atom* for $\{J_k\}_{k \geq 0}$ and then it is a classical result that the returns of $\{J_k\}_{k \geq 0}$ to 0 are not only $\mathcal{P}^a_\varepsilon$-a.s. finite, but also integrable:

$$\mathcal{P}^a_\varepsilon(\zeta_1 < \infty) = 1, \qquad \mathcal{E}^a_\varepsilon(\zeta_1) = \frac{1}{\nu^a_\varepsilon(\{0\})} < \infty, \qquad (\varepsilon \geq \varepsilon^a_c); \tag{5.5}$$

cf. [26], Chapter 5. Therefore also $\chi_1 = \tau_{\zeta_1}$ is $\mathcal{P}^a_\varepsilon$-a.s. a finite random variable for $\varepsilon \geq \varepsilon^a_c$. □

5.2. *Proof of Theorem 1.2.* We start showing that the quantities $\mathrm{F}^a(\varepsilon)$ and $\varepsilon^a_c$, that were defined in (4.10) and appear in the definition (3.9) of the kernel $\mathsf{K}^{a,\varepsilon}_{x,dy}(n)$, indeed coincide with the corresponding quantities defined in the Introduction. By (5.3) we have

$$\frac{1}{N} \log \mathcal{Z}^a_{\varepsilon,N} = \mathrm{F}^a(\varepsilon) - \frac{2}{N} \log \varepsilon + \frac{1}{N} \log \mathcal{P}^a_\varepsilon(N + 1 \in \chi).$$

Since by (1.9) we have $\liminf_{N\to\infty} N^{-1}\log \mathcal{Z}^a_{\varepsilon,N} \geq 0$, when $\mathrm{F}^a(\varepsilon)=0$ the trivial bound $\mathcal{P}^a_\varepsilon(N+1\in\chi)\leq 1$ shows that (1.8) holds true. Therefore to complete the identification of $\mathrm{F}^a$ it suffices to show that when $\mathrm{F}^a(\varepsilon)>0$, that is, when $\varepsilon>\varepsilon^a_c$, we have

$$\lim_{N\to\infty}\frac{1}{N}\log \mathcal{P}^a_\varepsilon(N+1\in\chi)=0. \tag{5.6}$$

However, it is well known (and easy to prove) that this relation is true in complete generality for any nonterminating aperiodic renewal process, and Proposition 5.1 shows that for $\varepsilon>\varepsilon^a_c$ the process $\{\chi_k\}_k$ under $\mathcal{P}^a_\varepsilon$ is indeed a genuine nonterminating renewal process, which is aperiodic because $\mathcal{P}^a_\varepsilon(\chi_1=1)>0$. This completes the identification of $\mathrm{F}^a(\varepsilon)$, as defined in (4.10), with the free energy defined by (1.8).

Completing the proof of Theorem 1.2 is now easy. By definition we have $B^{\mathrm{p},0}_{x,dy}\geq B^{\mathrm{w},0}_{x,dy}$ and one checks easily that the nonnegative kernel $B^{\mathrm{p},0}_{x,dy}-B^{\mathrm{w},0}_{x,dy}$ has strictly positive spectral radius. Then Lemma 4.2 with $A_{x,dy}=B^{\mathrm{w},0}_{x,dy}$ and $C_{x,dy}=B^{\mathrm{p},0}_{x,dy}-B^{\mathrm{w},0}_{x,dy}$ yields $\delta^{\mathrm{w}}(0)<\delta^{\mathrm{p}}(0)$, that is, $\varepsilon^{\mathrm{w}}_c>\varepsilon^{\mathrm{p}}_c$ by (4.10). The analyticity of $\mathrm{F}^a(\cdot)$ on $(\varepsilon^a_c,\infty)$ has been already discussed in Section 4.2. It remains to prove (1.11). Note that $e^{-\lambda}\mathsf{F}^a_{x,dy}(1)\leq B^{a,\lambda}_{x,dy}\leq e^{-\lambda}B^{a,0}_{x,dy}$ by (4.4), hence $c\,e^{-\lambda}\leq\delta^a(\lambda)\leq c'\,e^{-\lambda}$ by (4.9), with $c,c'>0$. Taking $\lambda=(\delta^a)^{-1}(1/\varepsilon)$ and using (4.10) we finally obtain $\mathrm{F}^a(\varepsilon)\sim\log(\varepsilon)$ as $\varepsilon\to\infty$.

## 6. Proof of Theorem 1.3.

In this section we prove Theorem 1.3, that is, we show that (1.16) holds true:

$$\liminf_{N\to\infty}\mathbb{E}^{\mathrm{w}}_{\varepsilon^{\mathrm{w}}_c,N}\left(\frac{\ell_N}{N}\right)>0,$$

where we have used (1.13). We introduce the quantity

$$\iota_N:=\#\{\chi\cap[0,N]\}=\max\{k\geq 0:\chi_k\leq N\}, \tag{6.1}$$

and since $\ell_N\geq\iota_N$, by Proposition 3.1 it is sufficient to show that

$$\liminf_{N\to\infty}\mathbb{E}^{\mathrm{w}}_{\varepsilon^{\mathrm{w}}_c,N}\left(\frac{\iota_N}{N}\right)=\liminf_{N\to\infty}\mathcal{E}^{\mathrm{w}}_{\varepsilon^{\mathrm{w}}_c}\left(\frac{\iota_N}{N}\,\middle|\,\mathcal{A}_N\right)>0. \tag{6.2}$$

We recall that the process $\{\chi_k\}_k$ under $\mathcal{P}^{\mathrm{w}}_{\varepsilon^{\mathrm{w}}_c}$ is a classical aperiodic renewal process by Proposition 5.1. Moreover we claim that $b:=\mathcal{E}^{\mathrm{w}}_{\varepsilon^{\mathrm{w}}_c}(\chi_1)<\infty$. Then by the strong law of large numbers we have $\iota_N/N\to 1/b$, $\mathcal{P}^{\mathrm{w}}_{\varepsilon^{\mathrm{w}}_c}$-a.s., and by the renewal theorem $\mathcal{P}^{\mathrm{w}}_{\varepsilon^{\mathrm{w}}_c}(\mathcal{A}_N)=\mathcal{P}^{\mathrm{w}}_{\varepsilon^{\mathrm{w}}_c}(N+1\in\chi)\to 1/b>0$ as $N\to\infty$. It follows that

$$\mathcal{E}^{\mathrm{w}}_{\varepsilon^{\mathrm{w}}_c}\left(\frac{\iota_N}{N}\,\middle|\,\mathcal{A}_N\right)\longrightarrow\frac{1}{b}>0\qquad(N\to\infty),$$

and (6.2) follows. It only remains to check that $\mathcal{E}^{\mathrm{w}}_{\varepsilon^{\mathrm{w}}_c}(\chi_1)<\infty$.

6.1. *A formula for* $\mathcal{E}^{\mathrm{w}}_{\varepsilon_c^{\mathrm{w}}}(\chi_1)$. The dependency on $\varepsilon_c^{\mathrm{w}}$ will be omitted from now on for notational convenience. We recall that $\chi_1 := \tau_1 + \cdots + \tau_{\zeta_1}$, where $\zeta_1 := \inf\{n \geq 0 \colon J_n = 0\}$. We introduce the kernel $\widehat{\mathsf{K}}^{\mathrm{w}}_{x,dy}(n) := \mathsf{K}^{\mathrm{w}}_{x,dy}(n)\mathbf{1}_{(y\neq 0)} (= \mathsf{K}^{\mathrm{w}}_{x,dy}(n)\,\mathbf{1}_{(n\geq 2)})$ that gives the transition probabilities of the process $\{(\tau_k, J_k)\}_k$ before the chain $\{J_k\}_k$ comes back to zero. Summing over the possible values of the variable $\zeta_1$, we obtain the expression

$$\mathcal{P}^{\mathrm{w}}(\chi_1 = n) = \int_{y\in\mathbb{R}} \sum_{k=0}^{\infty} (\widehat{\mathsf{K}}^{\mathrm{w}})^{*k}_{0,dy}(n-1) \cdot \mathsf{K}^{\mathrm{w}}_{y,\{0\}}(1), \tag{6.3}$$

where the convolution $*$ of kernels is defined in Section 1.7. Now observe that

$$\begin{aligned}
&\sum_{m\in\mathbb{N}} m \cdot \sum_{k=0}^{\infty} (\widehat{\mathsf{K}}^{\mathrm{w}})^{*k}_{0,dy}(m) \\
&\quad = \sum_{k=1}^{\infty}\sum_{m=1}^{\infty} \sum_{\substack{t_1,\ldots,t_k\in\mathbb{N} \\ t_1+\cdots+t_k=m}} (t_1+\cdots+t_k)(\widehat{\mathsf{K}}^{\mathrm{w}}(t_1)\circ\cdots\circ\widehat{\mathsf{K}}^{\mathrm{w}}(t_k))_{0,\,dy} \\
&\quad = \sum_{k=1}^{\infty}\sum_{i=1}^{k} ((\widehat{D}^{\mathrm{w}})^{\circ(i-1)} \circ \widehat{M}^{\mathrm{w}} \circ (\widehat{D}^{\mathrm{w}})^{\circ(k-i)})_{0,\,dy} \\
&\quad = ((1-\widehat{D}^{\mathrm{w}})^{-1} \circ \widehat{M}^{\mathrm{w}} \circ (1-\widehat{D}^{\mathrm{w}})^{-1})_{0,\,dy},
\end{aligned}$$

where $\widehat{D}^{\mathrm{w}}_{x,dy} := \sum_{n\in\mathbb{N}} \widehat{\mathsf{K}}^{\mathrm{w}}_{x,dy}(n) = D^{\mathrm{w}}_{x,dy}\mathbf{1}_{(y\neq 0)}$ and $\widehat{M}^{\mathrm{w}}_{x,dy} := \sum_{n\in\mathbb{N}} n \cdot \widehat{\mathsf{K}}^{\mathrm{w}}_{x,dy}(n)$. Notice that $\mathsf{K}^{\mathrm{w}}_{y,\{0\}}(1) = D^{\mathrm{w}}_{y,\{0\}}$ by (3.13), therefore by (6.3) we have

$$\begin{aligned}
\mathcal{E}^{\mathrm{w}}(\chi_1) &= \sum_{n\in\mathbb{N}} n \cdot \mathcal{P}^{\mathrm{w}}(\chi_1 = n) = 1 + \sum_{n\in\mathbb{N}} (n-1)\cdot\mathcal{P}^{\mathrm{w}}(\chi_1 = n) \\
&= 1 + ((1-\widehat{D}^{\mathrm{w}})^{-1} \circ \widehat{M}^{\mathrm{w}} \circ (1-\widehat{D}^{\mathrm{w}})^{-1} \circ D^{\mathrm{w}})_{0,\{0\}}.
\end{aligned}$$

We recall that $D^{\mathrm{w}}_{x,dy}$ is the transition kernel of the Markov chain $\{J_k\}_k$ under $\mathcal{P}^{\mathrm{w}}$, which is positive recurrent with invariant probability measure $\nu^{\mathrm{w}}(\cdot) = \nu_{\varepsilon_c^{\mathrm{w}}}^{\mathrm{w}}(\cdot)$ defined in (4.13). Since $\nu^{\mathrm{w}}(\{0\}) > 0$, the state 0 is an atom for $\{J_k\}_k$ and therefore by [26] the following formulas hold: for all $x, y \in \mathbb{R}$

$$\begin{aligned}
(1-\widehat{D}^{\mathrm{w}})^{-1}_{0,\,dx} &= \frac{\nu^{\mathrm{w}}(dx)}{\nu^{\mathrm{w}}(\{0\})} = \frac{v^{\mathrm{w}}(x)w^{\mathrm{w}}(x)}{v^{\mathrm{w}}(0)w^{\mathrm{w}}(0)}\mu(dx), \\
((1-\widehat{D}^{\mathrm{w}})^{-1} \circ D^{\mathrm{w}})_{y,\{0\}} &= 1.
\end{aligned}$$

Then we finally come to the expression

$$\mathcal{E}^{\mathrm{w}}(\chi_1) = 1 + \varepsilon_c^{\mathrm{w}} \int_{x,y\in\mathbb{R}^+} \mu(dx)\,\mu(dy)\,\frac{w^{\mathrm{w}}(x)}{w^{\mathrm{w}}(0)} \left(\sum_{n\in\mathbb{N}} n\,\mathsf{f}^{\mathrm{w}}_{x,y}(n)\right) \frac{v^{\mathrm{w}}(y)}{v^{\mathrm{w}}(0)}. \tag{6.4}$$

6.2. *Finiteness of* $\mathcal{E}_{\varepsilon_c^{\mathrm{w}}}(\chi_1)$. First we state two relations that we prove below: for both $a \in \{\mathrm{p}, \mathrm{w}\}$ and for all $\varepsilon > 0$:

$$\begin{aligned} &\int_{\mathbb{R}} \mu(dx) v_\varepsilon^a(x) \, \log(x^2 \vee 1) < \infty, \\ &\int_{\mathbb{R}} \mu(dx) w_\varepsilon^a(x) \, \log(x^2 \vee 1) < \infty. \end{aligned} \tag{6.5}$$

We aim at showing that the r.h.s. of (6.4) is finite. We start considering the terms in the sum with $n \geq (x^2 \vee y^2)$: applying (3.8), (2.12) and the first relation in (4.2) we obtain

$$\sum_{n \geq (x^2 \vee y^2)} n \mathsf{f}_{x,y}^{\mathrm{w}}(n) \leq \sum_{n \geq (x^2 \vee y^2)} \frac{(const.)}{n \, (\log n)^{c_+}} \leq (const.'),$$

because $c_+ > 1$ and hence the sum converges. Therefore the contribution to the r.h.s. of (6.4) is bounded by

$$\begin{aligned} &\varepsilon_c^{\mathrm{w}} (const.') \int_{x,y \in \mathbb{R}^+} \mu(dx)\mu(dy) \frac{w^{\mathrm{w}}(x)}{w^{\mathrm{w}}(0)} \frac{v^{\mathrm{w}}(y)}{v^{\mathrm{w}}(0)} \\ &\qquad \leq \varepsilon_c^{\mathrm{w}} (const.') \frac{\|v^{\mathrm{w}}\|_1 \|w^{\mathrm{w}}\|_1}{v^{\mathrm{w}}(0) w^{\mathrm{w}}(0)}, \end{aligned}$$

where $\|\cdot\|_1$ denotes the norm in $L^1(\mathbb{R}, d\mu)$. Notice that $\|v^{\mathrm{w}}\|_1 < \infty$, $\|w^{\mathrm{w}}\|_1 < \infty$ by (6.5).

Next we deal with the terms in the r.h.s. of (6.4) with $n < (x^2 \vee y^2)$. From the bound $\mathsf{f}_{x,y}^{\mathrm{w}}(n) \leq \mathsf{f}_{x,y}^{\mathrm{p}}(n) \leq C/n^2$ [cf. (4.2)], we have

$$\sum_{n < (x^2 \vee y^2)} n \mathsf{f}_{x,y}^{\mathrm{w}}(n) \leq C \sum_{n < (x^2 \vee y^2)} \frac{1}{n} \leq C(\log(x^2 \vee 1) + \log(y^2 \vee 1)),$$

and using again (6.5) we see that the r.h.s. of (6.4) is indeed finite.

6.3. *Proof of (6.5).* We focus on the first relation, the second one being analogous. By (4.12) we have for both $a \in \{\mathrm{p}, \mathrm{w}\}$

$$v_\varepsilon^a(x) \leq (const.) \sum_{n \in \mathbb{N}} \int_{y \in \mathbb{R}} \mathsf{f}_{x,y}^{\mathrm{p}}(n) \, v_\varepsilon^a(y) \mu(dy),$$

because $\mathsf{f}_{x,y}^{\mathrm{w}}(n) \leq \mathsf{f}_{x,y}^{\mathrm{p}}(n)$. Setting $u_n(y) := \int_{x \in \mathbb{R}} \mathsf{f}_{x,y}^{\mathrm{p}}(n) \log(x^2 \vee 1) \mu(dx)$ and applying the Cauchy–Schwarz inequality we get

$$\begin{aligned} \int_{x \in \mathbb{R}} \mu(dx) v_\varepsilon^a(x) \log(x^2 \vee 1) &\leq (const.) \sum_{n \in \mathbb{N}} \int_{y \in \mathbb{R}} \mu(dy) \, u_n(y) v_\varepsilon^a(y) \\ &\leq (const.) \|v_\varepsilon^a\|_2 \cdot \left( \sum_{n \in \mathbb{N}} \|u_n\|_2 \right), \end{aligned} \tag{6.6}$$

where $\|v_\varepsilon^a\|_2 := (\int_{\mathbb{R}} (v_\varepsilon^a(x))^2 \mu(dx))^{1/2}$ and likewise for $\|u_n\|_2$. Setting for short $\mathsf{f}^{\mathrm{p}}_{\mathbb{R},y}(n) := \int_{\mathbb{R}} \mathsf{f}^{\mathrm{p}}_{x,y}(n)\, \mu(dx)$, by Jensen's inequality we have

$$\begin{aligned}\|u_n\|_2^2 = \int_{y\in\mathbb{R}} (u_n(y))^2 \mu(dy) &= \int_{y\in\mathbb{R}} \mu(dy) \left( \int_{x\in\mathbb{R}} \mathsf{f}^{\mathrm{p}}_{x,y}(n) \log(x^2 \vee 1) \mu(dx) \right)^2 \\ &\le \int_{y\in\mathbb{R}} \mu(dy) (\mathsf{f}^{\mathrm{p}}_{\mathbb{R},y}(n))^2 \int_{x\in\mathbb{R}} \frac{\mathsf{f}^{\mathrm{p}}_{x,y}(n)}{\mathsf{f}^{\mathrm{p}}_{\mathbb{R},y}(n)} \log^2(x^2 \vee 1) \mu(dx),\end{aligned}$$

and since $\mathsf{f}^{\mathrm{p}}_{\mathbb{R},y}(n) \le C/n^{3/2}$ by (4.2), Fubini's theorem yields

$$\int_{y\in\mathbb{R}} (u_n(y))^2 \mu(dy) \le \frac{C}{n^{3/2}} \frac{1}{n} \int_{x\in\mathbb{R}} n \cdot \mathsf{f}^{\mathrm{p}}_{x,\mathbb{R}}(n) \log^2(x^2 \vee 1) \mu(dx).$$

Observe that $\mathsf{f}^{\mathrm{p}}_{\mathbb{R},\mathbb{R}}(n) := \int_{x\in\mathbb{R}} \mathsf{f}^{\mathrm{p}}_{x,\mathbb{R}}(n) \mu(dx) = 1/n$ by (4.3), therefore $x \mapsto n \cdot \mathsf{f}^{\mathrm{p}}_{x,\mathbb{R}}(n)$ is a probability density; in fact it is the density of the random variable $Z_n/n$ under $\mathbf{P}^{(0,0)}$, as follows from (3.7). Since $\log^2(x^2 \vee 1) \le \log^2(x^2 \vee e)$ and the function $z \mapsto \log^2(z)$ is concave on the half-line $[e, \infty)$, by Jensen's inequality we get

$$\begin{aligned}\int_{x\in\mathbb{R}} n \cdot \mathsf{f}^{\mathrm{p}}_{x,\mathbb{R}}(n) \log^2(x^2 \vee 1) \mu(dx) &\le \log^2 \left( \int_{x\in\mathbb{R}} n \cdot \mathsf{f}^{\mathrm{p}}_{x,\mathbb{R}}(n) (x^2 \vee e)\, \mu(dx) \right) \\ &\le \log^2 \left( \mathbf{E}^{(0,0)} \left( \frac{Z_n^2}{n^2} \vee e \right) \right) \\ &\le \log^2 (e + (const.) \cdot n),\end{aligned}$$

because by (2.2) we have $\mathbf{E}^{(0,0)}(Z_n^2) \sim \sigma^2 \cdot n^3/3$ as $n \to \infty$. It follows that

$$\int_{y\in\mathbb{R}} (u_n(y))^2 \mu(dy) \le \frac{C}{n^{5/2}} \log^2(e + (const.) \cdot n)$$

and therefore

$$\sum_{n\in\mathbb{N}} \|u_n\|_2 = \sum_{n\in\mathbb{N}} \sqrt{\int_{y\in\mathbb{R}} (u_n(y))^2 \mu(dy)} < \infty.$$

Looking back to (6.6), equation (6.5) is proven.

**7. Proof of Theorem 1.4.** In this section we prove Theorem 1.4. We start proving (1.18), that we rewrite for convenience: there exists a positive constant $c_5$ such that for small $\delta > 0$

$$\textsc{f}^{\mathrm{p}}(\varepsilon_c^{\mathrm{p}} + \delta) \ge c_5 \frac{\delta}{\log 1/\delta}. \tag{7.1}$$

To this purpose, we observe that the following relation holds true for all $\varepsilon, \varepsilon' > 0$:

$$\frac{\mathcal{Z}^{\mathrm{p}}_{\varepsilon',N}}{\mathcal{Z}^{\mathrm{p}}_{\varepsilon,N}} = \mathbb{E}^{\mathrm{p}}_{\varepsilon,N}\left(\left(\frac{\varepsilon'}{\varepsilon}\right)^{\ell_N - 1}\right),$$

where we recall that $\ell_N = \#\{i \in \{1,\dots,N\} : \varphi_i = 0\}$. This relation follows easily from the definition (1.5) of our model, because the differences between $\mathcal{Z}^{\mathrm{p}}_{\varepsilon,N}$ and $\mathcal{Z}^{\mathrm{p}}_{\varepsilon',N}$ are only given by the different size of the pinning reward for each time the path touches the $x$-axis (see also Appendix A). Then we have the easy lower bound

$$\frac{\mathcal{Z}^{\mathrm{p}}_{\varepsilon',N}}{\mathcal{Z}^{\mathrm{p}}_{\varepsilon,N}} \ge \mathbb{E}^{\mathrm{p}}_{\varepsilon,N}\left(\left(\frac{\varepsilon'}{\varepsilon}\right)^{\iota_N - 1}\right) = \mathcal{E}^{\mathrm{p}}_{\varepsilon}\left(\left(\frac{\varepsilon'}{\varepsilon}\right)^{\iota_N - 1} \,\middle|\, N+1 \in \chi\right),$$

where we recall that $\iota_N = \#\{i \in \{1,\dots,N\} : \chi_i = 0\}$ and we have applied Proposition 3.1. Setting $\varepsilon = \varepsilon_c^{\mathrm{p}}$, $\varepsilon' = \varepsilon_c^{\mathrm{p}} + \delta$, taking logarithms and letting $N \to \infty$ yields

$$\begin{aligned} \mathrm{F}^{\mathrm{p}}(\varepsilon_c^{\mathrm{p}} + \delta) \;&\ge\; \mathrm{F}^{\mathrm{p}}(\varepsilon_c^{\mathrm{p}}) + \limsup_{N\to\infty} \frac{1}{N} \log \mathcal{E}^{\mathrm{p}}_{\varepsilon_c^{p}}\left(\left(\frac{\varepsilon_c^{\mathrm{p}} + \delta}{\varepsilon_c^{\mathrm{p}}}\right)^{\iota_N - 1} \,\middle|\, N+1 \in \chi\right) \\ &=: \mathrm{F}^{\mathrm{p}}(\varepsilon_c^{\mathrm{p}}) + \mathrm{G}(\delta), \end{aligned}$$

and since $\mathrm{F}^{\mathrm{p}}(\varepsilon_c^{\mathrm{p}}) = 0$ we get

$$\mathrm{F}^{\mathrm{p}}(\varepsilon_c^{\mathrm{p}} + \delta) \ge \mathrm{G}(\delta).$$

Now notice that under $\mathcal{P}^{\mathrm{p}}_{\varepsilon_c^{p}}$ the process $\{\chi_k\}_k$ is a classical, nonterminating renewal process, by Proposition 5.1. Therefore $\mathrm{G}(\delta)$ is just the free energy of a *classical pinning model*, that is, a renewal process rewarded $\frac{\varepsilon_c^{\mathrm{p}}+\delta}{\varepsilon_c^{\mathrm{p}}}$ at each renewal epoch. The asymptotic behavior of $\mathrm{G}(\delta)$ is then given by Theorem 2.1 of [13]: since $\mathcal{P}^{\mathrm{p}}_{\varepsilon_c^{p}}(\chi_1 = n) \sim c/n^2$ as $n \to \infty$, as we prove in Proposition 7.1 below, we have $\mathrm{G}(\delta) \sim c'\delta/\log(1/\delta)$ as $\delta \downarrow 0$, and (7.1) is proven.

It only remains to show that (1.17) holds true, that is,

$$\limsup_{\varepsilon \downarrow \varepsilon_c^{\mathrm{p}}} \limsup_{N\to\infty} \mathcal{E}^{\mathrm{p}}_{\varepsilon}\left(\frac{\ell_N}{N} \,\middle|\, \mathcal{A}_N\right) = 0, \tag{7.2}$$

where we have used (1.13) and Proposition 3.1. The idea is to focus first on the variable $\iota_N$, which is easier to handle, and then to make the comparison with $\ell_N$.

7.1. *From $\iota_N$ to $\ell_N$.* We recall that the process $\{\chi_k\}_k$ under the law $\mathcal{P}^{\mathrm{p}}_{\varepsilon}$ is a classical aperiodic renewal process; see Proposition 5.1. We introduce the step law

$$q_\varepsilon(n) := \mathcal{P}^{\mathrm{p}}_{\varepsilon}(\chi_1 = n), \tag{7.3}$$

whose asymptotic behavior as $n \to \infty$, when $\varepsilon$ is close to $\varepsilon_c^{\mathrm{p}}$, is given by the following:

PROPOSITION 7.1. *There exists $\alpha > 0$ such that for every $\varepsilon \in [\varepsilon_c^{\mathrm{p}}, \varepsilon_c^{\mathrm{p}} + \alpha]$ we have*

$$q_\varepsilon(n) \sim \frac{C_\varepsilon}{n^2} \exp(-\mathrm{F}^{\mathrm{p}}(\varepsilon) \cdot n) \qquad (n \to \infty), \tag{7.4}$$

*where $C_\varepsilon \in (0, \infty)$ is a continuous function of $\varepsilon$.*

We postpone the proof to the next paragraphs; for the moment we focus on the consequences. We assume in the following that $\varepsilon_c^{\mathrm{p}} < \varepsilon \le \varepsilon_c^{\mathrm{p}} + \alpha$.

Let us set $G_\varepsilon := \mathcal{E}_\varepsilon^{\mathrm{p}}(\chi_1) < \infty$ by Proposition 7.1. A standard Tauberian theorem (cf. [5], Theorem 1.7.1) gives the asymptotic behavior of $G_\varepsilon$ as $\varepsilon \downarrow \varepsilon_c^{\mathrm{p}}$:

$$G_\varepsilon \sim (const.) \log \frac{1}{\mathrm{F}^{\mathrm{p}}(\varepsilon)} \qquad (\varepsilon \downarrow \varepsilon_c^{\mathrm{p}}). \tag{7.5}$$

Notice that the classical Renewal theorem yields

$$\mathcal{P}_\varepsilon^{\mathrm{p}}(\mathcal{A}_N) = \mathcal{P}_\varepsilon^{\mathrm{p}}(N+1 \in \chi) \to \frac{1}{G_\varepsilon} > 0 \qquad (N \to \infty). \tag{7.6}$$

Therefore by the weak law of large numbers for the process $\{\chi_k\}_k$ we have

$$\begin{aligned} \mathcal{P}_\varepsilon^{\mathrm{p}}\left(\frac{\iota_N}{N} \ge \frac{2}{G_\varepsilon} \,\Big|\, \mathcal{A}_N\right) &= \frac{\mathcal{P}_\varepsilon^{\mathrm{p}}(\chi_{\lfloor 2N/G_\varepsilon \rfloor} \le N)}{\mathcal{P}_\varepsilon^{\mathrm{p}}(\mathcal{A}_N)} \\ &= \frac{\mathcal{P}_\varepsilon^{\mathrm{p}}\ (\chi_{\lfloor 2N/G_\varepsilon \rfloor}/\lfloor 2N/G_\varepsilon \rfloor \le G_\varepsilon/2)}{\mathcal{P}_\varepsilon^{\mathrm{p}}(\mathcal{A}_N)} \stackrel{N\to\infty}{\longrightarrow} 0. \end{aligned} \tag{7.7}$$

We recall that $\zeta_k$ denotes the epoch of the $k$th return of the process $\{J_k\}_k$ to the state zero [cf. 5.4], and that $\{\zeta_k\}_{k \ge 0}$ under $\mathcal{P}_\varepsilon^{\mathrm{p}}$ is a nonterminating renewal process with finite mean $m_\varepsilon := \mathcal{E}_\varepsilon^{\mathrm{p}}(\zeta_1) = 1/\nu_\varepsilon^{\mathrm{p}}(\{0\}) < \infty$; cf. (5.5). Then by the weak law of large numbers

$$\mathcal{P}_\varepsilon^{\mathrm{p}}\left(\frac{\zeta_k}{k} > 2m_\varepsilon\right) \longrightarrow 0 \qquad (k \to \infty). \tag{7.8}$$

We are finally ready to estimate $\ell_N$. A trivial bound yields

$$\begin{aligned} \mathcal{P}_\varepsilon^{\mathrm{p}}\left(\frac{\ell_N}{N} \ge \frac{4m_\varepsilon}{G_\varepsilon} \,\Big|\, \mathcal{A}_N\right) &\le \mathcal{P}_\varepsilon^{\mathrm{p}}\left(\frac{\iota_N}{N} \ge \frac{2}{G_\varepsilon} \,\Big|\, \mathcal{A}_N\right) \\ &\quad + \mathcal{P}_\varepsilon^{\mathrm{p}}\left(\frac{\ell_N}{N} \ge \frac{4m_\varepsilon}{G_\varepsilon}, \frac{\iota_N}{N} < \frac{2}{G_\varepsilon} \,\Big|\, \mathcal{A}_N\right). \end{aligned}$$

The first term in the r.h.s. vanishes as $N \to \infty$ by (7.7). For the second term we observe that by definition $\ell_N = \zeta_{\iota_N}$ on the event $\mathcal{A}_N$, hence by an inclusion argument we have

$$\mathcal{P}^{\mathrm{p}}_{\varepsilon}\left( \frac{\ell_N}{N} \geq \frac{4m_\varepsilon}{G_\varepsilon}, \frac{\iota_N}{N} < \frac{2}{G_\varepsilon} \,\Big|\, \mathcal{A}_N \right) \leq \mathcal{P}^{\mathrm{p}}_{\varepsilon}\left( \frac{\zeta_{\lfloor 2N/G_\varepsilon \rfloor}}{\lfloor 2N/G_\varepsilon \rfloor} > 2m_\varepsilon \,\Big|\, \mathcal{A}_N \right) \stackrel{N\to\infty}{\longrightarrow} 0\,,$$

having used (7.8) and (7.6). Since $\ell_N/N \leq 1$, we can finally write

$$\limsup_{N\to\infty} \mathcal{E}^{\mathrm{p}}_{\varepsilon}\left( \frac{\ell_N}{N} \,\Big|\, \mathcal{A}_N \right) \leq \frac{4m_\varepsilon}{G_\varepsilon} + \limsup_{N\to\infty} \mathcal{P}^{\mathrm{p}}_{\varepsilon}\left( \frac{\ell_N}{N} \geq \frac{4m_\varepsilon}{G_\varepsilon} \,\Big|\, \mathcal{A}_N \right) = \frac{4m_\varepsilon}{G_\varepsilon}.$$

Now observe that as $\varepsilon \downarrow \varepsilon_c^{\mathrm{p}}$ we have $m_\varepsilon \to 1/\nu_{\varepsilon_c^{\mathrm{p}}}(\{0\}) < \infty$ and moreover $G_\varepsilon \to +\infty$ by (7.5). Then we let $\varepsilon \downarrow \varepsilon_c^{\mathrm{p}}$ in the last equation and (7.2) is proven.

7.2. *A Markov renewal theorem with infinite mean.* Before proving Proposition 7.1, we derive a generalized renewal theorem in our Markovian setting. Since the steps are more transparent if carried out in a general setting, we assume that we are given a kernel $\mathsf{A}_{x,dy}(n)$ satisfying the following assumptions:

1. the spectral radius of $G_{x,dy} := \sum_{n\in\mathbb{N}} \mathsf{A}_{x,dy}(n)$ is strictly less than 1;
2. we have $\mathsf{A}_{x,dy}(n) \sim L_{x,dy}/n^2$ as $n \to \infty$, for some kernel $L_{x,dy}$ (for the precise meaning of this relation we refer to Section 1.7), and moreover $\mathsf{A}_{x,dy}(n) \leq cL_{x,dy}/n^2$ for every $n \in \mathbb{N}$, where $c$ is a positive constant;
3. there exists $\beta > 1$ such that $((1-\beta G)^{-1} \circ L \circ (1-\beta G)^{-1})_{x,F} < \infty$ for every $x \in \mathbb{R}$ and for every *bounded* Borel set $F \subset \mathbb{R}$.

The result we are going to prove is the following asymptotic relation:

$$\sum_{k=0}^{\infty} \mathsf{A}^{*k}_{x,dy}(n) \sim \frac{((1-G)^{-1} \circ L \circ (1-G)^{-1})_{x,dy}}{n^2} \qquad (n \to \infty). \tag{7.9}$$

The path we follow is close to [8], Section 3.4. We start proving by induction the following bound: for all $k, n \in \mathbb{N}$ and $x, y \in \mathbb{R}$

$$\mathsf{A}^{*k}_{x,dy}(n) \leq c\, k^2\, \frac{\sum_{i=0}^{k-1} (G^{\circ i} \circ L \circ G^{\circ[(k-1)-i]})_{x,dy}}{n^2}. \tag{7.10}$$

The $k = 1$ case holds by assumption (2). Then we consider the even-$k$ case: by the definition of the convolution $*$ we have

$$\mathsf{A}^{*(2k)}_{x,dy}(n) \leq \sum_{h=1}^{\lceil n/2 \rceil} \int_{z\in\mathbb{R}} (\mathsf{A}^{*k}_{x,dz}(h) \cdot \mathsf{A}^{*k}_{z,dy}(n-h) + \mathsf{A}^{*k}_{x,dz}(n-h) \cdot \mathsf{A}^{*k}_{z,dy}(h)).$$

Observing that $\sum_{h\in\mathbb{N}} \mathsf{A}^{*k}_{x,dy}(h) = G^{\circ k}_{x,dy}$ and applying the inductive step we get

$$\begin{aligned}\mathsf{A}^{*(2k)}_{x,dy}(n) &\le \frac{c\,k^2}{(n/2)^2}\Bigg[\Bigg(G^{\circ k}\circ\Bigg(\sum_{i=0}^{k-1} G^{\circ i}\circ L\circ G^{\circ[(k-1)-i]}\Bigg)\Bigg)_{x,dy} \\ &\qquad + \Bigg(\Bigg(\sum_{i=0}^{k-1} G^{\circ i}\circ L\circ G^{\circ[(k-1)-i]}\Bigg)\circ G^{\circ k}\Bigg)_{x,dy}\Bigg] \\ &= \frac{c(2k)^2}{n^2}\sum_{i=0}^{2k-1}(G^{\circ i}\circ L\circ G^{\circ[(2k-1)-i]})_{x,dy},\end{aligned}$$

so that (7.10) is proven, the odd-$k$ case being analogous. Note that, choosing a constant $c'>0$ such that $k^2\le c'\beta^k$ for every $k$, assumption (3) yields that for every $x\in\mathbb{R}$ and for every bounded Borel set $F\subset\mathbb{R}$

$$\begin{aligned}&\sum_{k=1}^{\infty}\sum_{i=0}^{k-1} k^2 (G^{\circ i}\circ L\circ G^{\circ[(k-1)-i]})_{x,F} \\ &\qquad\le c'\sum_{k=1}^{\infty}\sum_{i=0}^{k-1}((\beta G)^{\circ i}\circ L\circ(\beta G)^{\circ[(k-1)-i]})_{x,F} \\ &\qquad = c'((1-\beta G)^{-1}\circ L\circ(1-\beta G)^{-1})_{x,F}<\infty.\end{aligned} \tag{7.11}$$

Next we claim that

$$\mathsf{A}^{*k}_{x,dy}(n)\sim\frac{\sum_{i=0}^{k-1}(G^{\circ i}\circ L\circ G^{\circ[(k-1)-i]})_{x,dy}}{n^2}\qquad (n\to\infty). \tag{7.12}$$

We proceed by induction: the $k=1$ case holds by assumption (1), while for general $k$

$$\begin{aligned}\mathsf{A}^{*k}_{x,dy}(n) = \sum_{h=1}^{n/2}\int_{z\in\mathbb{R}}(&\mathsf{A}^{*(k-1)}_{x,dz}(h)\cdot\mathsf{A}_{z,dy}(n-h) \\ &+\mathsf{A}^{*(k-1)}_{x,dz}(n-h)\cdot\mathsf{A}_{z,dy}(h)).\end{aligned}$$

Applying the induction step and using dominated convergence, thanks to (7.10) and (7.11), we have [observe that $\sum_{h\in\mathbb{N}}\mathsf{A}^{*m}_{x,dy}(h) = G^{\circ m}_{x,dy}$]

$$\begin{aligned}n^2\mathsf{A}^{*k}_{x,dy}(n) &\overset{n\to\infty}{\longrightarrow}\Bigg(G^{\circ(k-1)}\circ L+\Bigg(\sum_{i=0}^{k-2}G^{\circ i}\circ L\circ G^{\circ[(k-2)-i]}\Bigg)\circ G\Bigg)_{x,\,dy} \\ &= \sum_{i=0}^{k-1}(G^{\circ i}\circ L\circ G^{\circ[(k-1)-i]})_{x,\,dy},\end{aligned}$$

and (7.12) is proven. Finally we can write as $n \to \infty$

$$\begin{aligned} n^2 \sum_{k\geq 0} \mathsf{A}^{*k}_{x,dy}(n) &\to \sum_{k\geq 0} \sum_{i=0}^{k-1} (G^{\circ i} \circ L \circ G^{\circ[(k-1)-i]})_{x,dy} \\ &= ((1-G)^{-1} \circ L \circ (1-G)^{-1})_{x,dy}, \end{aligned}$$

where we have applied (7.12) and again dominated convergence, using (7.10) and (7.11). This completes the proof of (7.9).

7.3. *Proof of Proposition 7.1.* We start from a close analog of (6.3), namely

$$\begin{aligned} q_\varepsilon(n) &= \int_{y\in\mathbb{R}} \left( \sum_{k=0}^{\infty} (\widehat{\mathsf{K}}^{\mathrm{p},\varepsilon})^{*k}_{0,dy}(n-1) \right) \cdot \mathsf{K}^{\mathrm{p},\varepsilon}_{y,\{0\}}(1) \\ &= \varepsilon \cdot e^{-\mathrm{F}(\varepsilon)n} \cdot \int_{y\in\mathbb{R}} \left( \sum_{k=0}^{\infty} \varepsilon^k (\widehat{\mathsf{F}}^{\mathrm{p}})^{*k}_{0,dy}(n-1) \right) \cdot \mathsf{F}^{\mathrm{p}}_{y,\{0\}}(1), \end{aligned} \tag{7.13}$$

where we have set $\widehat{\mathsf{K}}^{\mathrm{p},\varepsilon}_{x,dy}(n) := \mathsf{K}^{\mathrm{p},\varepsilon}_{x,dy}(n)\mathbf{1}_{(y\neq 0)}$ and $\widehat{\mathsf{F}}^{\mathrm{p}}_{x,dy}(n) := \mathsf{F}^{\mathrm{p}}_{x,dy}(n)\mathbf{1}_{(y\neq 0)}$. What we need is the asymptotic behavior as $n \to \infty$ of the r.h.s. of (7.13) and to this purpose we are going to apply the results of Section 7.2 to the kernel $\mathsf{A}_{x,dy}(n) = \varepsilon \cdot \widehat{\mathsf{F}}^{\mathrm{p}}_{x,dy}(n)$.

We need to check that the assumptions (1)–(3) are satisfied. The asymptotic behavior of $\widehat{\mathsf{F}}^{\mathrm{p}}_{x,dy}(n)$ is obtained by (4.1):

$$\widehat{\mathsf{F}}^{\mathrm{p}}_{x,dy}(n) \sim \frac{c}{n^2}\, dy \qquad (n \to \infty), \tag{7.14}$$

and from (4.2) we see that assumption (2) is checked. We set for simplicity

$$\widehat{B}^{\mathrm{p}}_{x,dy} := \sum_{n\in\mathbb{N}} \widehat{F}^{\mathrm{p}}_{x,dy}(n) = B^{\mathrm{p},0}_{x,dy} \mathbf{1}_{(y\neq 0)}\,,$$

where the kernel $B^{\mathrm{p},0}_{x,dy}$ was defined in (4.4). The spectral radius of $\varepsilon_c^{\mathrm{p}} \cdot B^{\mathrm{p},0}_{x,dy}$ equals 1 by the very definition of $\varepsilon_c^{\mathrm{p}}$. Then applying Lemma 4.2 with $A_{x,dy} = \varepsilon_c^{\mathrm{p}} \cdot \widehat{B}^{\mathrm{p}}_{x,dy}$ and $C_{x,dy} = \varepsilon_c^{\mathrm{p}} \cdot B^{\mathrm{p},0}_{x,dy} - \varepsilon_c^{\mathrm{p}} \cdot \widehat{B}^{\mathrm{p}}_{x,dy}$ [it is easily seen that the spectral radius of $C_{x,dy} = \varepsilon_c^{\mathrm{p}} \cdot B^{\mathrm{p},0}_{x,dy}\mathbf{1}_{(y=0)} = \varepsilon_c^{\mathrm{p}} \cdot e^{-V(x)}\delta_0(dy)$ is strictly positive] we have that the spectral radius of $\varepsilon_c^{\mathrm{p}} \cdot \widehat{B}^{\mathrm{p}}_{x,dy}$ is strictly smaller than 1. By continuity there exists $\alpha > 0$ such that the spectral radius of $\varepsilon \cdot \widehat{B}^{\mathrm{p}}_{x,dy}$ is strictly smaller than 1 for every $\varepsilon \in [\varepsilon_c^{\mathrm{p}}, \varepsilon_c^{\mathrm{p}} + \alpha]$. Then assumption (1) is verified and it only remains to check assumption (3), that is,

$$\left( \int_{y\in\mathbb{R}} (1 - \beta\varepsilon \widehat{B}^{\mathrm{p}})^{-1}_{x,dy} \right) \left( \int_{z\in\mathbb{R}} \mu(dz) \int_{w\in\mathcal{A}} (1-\beta\varepsilon\widehat{B}^{\mathrm{p}})^{-1}_{z,dw}\, v(w) \right) < \infty, \tag{7.15}$$

for some $\beta > 1$. Let us focus on the first integral: we can write

$$\int_{y\in\mathbb{R}} (1-\beta\varepsilon\widehat{B}^{\mathrm{p}})^{-1}_{x,dy} = 1 + \sum_{n\geq 0} (\beta\varepsilon) \int_{z\in\mathbb{R}} (\beta\varepsilon\widehat{B}^{\mathrm{p}})^{\circ n}_{x,dz}\, g(z)\,,$$

where $g(z) := \int_{y\in\mathbb{R}} \widehat{B}^{\mathrm{p}}_{z,dy}$. We choose $\beta$ sufficiently close to 1 so that the spectral radius of $\beta\varepsilon\widehat{B}^{\mathrm{p}}$, let us call it $\rho$, is strictly smaller than 1. Denoting by $\|\cdot\|$ the operator norm in $L^2(\mathbb{R}, d\mu)$, a classical result gives $\|(\beta\varepsilon\widehat{B}^{\mathrm{p}})^{\circ n}\|^{1/n} \to \rho < 1$ as $n\to\infty$; cf. [20], Chapter III, Section 6.2. Therefore, if we show that $g(\cdot) \in L^2(\mathbb{R}, d\mu)$, we obtain

$$\int_{y\in\mathbb{R}} (1-\beta\varepsilon\widehat{B}^{\mathrm{p}})^{-1}_{x,dy} \leq 1 + (\beta\varepsilon) \sum_{n\geq 0} \|(\beta\varepsilon\widehat{B}^{\mathrm{p}})^{\circ n}\|\|g\| < \infty.$$

To prove that $g(\cdot) \in L^2(\mathbb{R}, d\mu)$, we observe that $g(z) = \sum_{n\in\mathbb{N}} \int_{y\in\mathbb{R}} \widehat{\mathsf{F}}^{\mathrm{p}}_{z,dy}(n)$ so that

$$\int_{x\in\mathbb{R}} g(x)^2 \mu(dx) = \sum_{n\in\mathbb{N}} \sum_{m\in\mathbb{N}} \int_{x\in\mathbb{R}} \biggl(\int_{y\in\mathbb{R}} \widehat{\mathsf{F}}^{\mathrm{p}}_{x,dy}(n)\biggr)\biggl(\int_{z\in\mathbb{R}} \widehat{\mathsf{F}}^{\mathrm{p}}_{x,dz}(m)\biggr)\mu(dx).$$

Using the symmetry in $n, m$ and applying (4.2) and (4.3), we finally obtain

$$\int_{x\in\mathbb{R}} g(x)^2 \mu(dx) \leq 2 \sum_{m\in\mathbb{N}} \sum_{n\geq m} \frac{C}{n^{3/2} m} < \infty.$$

With similar arguments one shows that also the second integral term in (7.15) is finite.

We can finally apply (7.9) in our setting, getting

$$\sum_{k=0}^{\infty} \varepsilon^k (\mathsf{F}^{\mathrm{p}})^{*k}_{0,dy}(n) \sim \frac{c\varepsilon}{n^2} \biggl(\int_{z\in\mathbb{R}} (1-\varepsilon\widehat{B}^{\mathrm{p}})^{-1}_{0,dz}\biggr) \cdot \biggl(\int_{x\in\mathbb{R}} dx\, (1-\varepsilon\widehat{B}^{\mathrm{p}})^{-1}_{x,dy}\biggr) \qquad (n\to\infty).$$

Coming back to (7.13), we can apply dominated convergence thanks to (7.10) and (7.11); we thus obtain (7.4), with $C_\varepsilon$ given by

$$C_\varepsilon := c\varepsilon^2 \biggl(\int_{z\in\mathbb{R}} (1-\varepsilon\widehat{B}^{\mathrm{p}})^{-1}_{0,dz}\biggr) \cdot \biggl(\int_{x,y\in\mathbb{R}} dx (1-\varepsilon\widehat{B}^{\mathrm{p}})^{-1}_{x,dy} e^{-V(y)}\biggr),$$

and the proof is complete.

## APPENDIX A: CONVEXITY OF THE FREE ENERGY

Recall the definition (1.12) of the contact number $\ell_N$ and observe that in any case $\ell_N \geq 1$ under $\mathbb{P}^a_{\varepsilon,N}$, because $\varphi_N = 0$. Setting $\Omega^{\mathrm{p}} := \mathbb{R}$ and $\Omega^{\mathrm{w}} := \mathbb{R}^+$, from the definitions (1.5) and (1.7) of our models we can write for $k \in \mathbb{Z}^+$

$$\begin{aligned}
&\mathbb{P}^a_{\varepsilon,N}(\ell_N = k+1) \\
&\quad = \frac{\varepsilon^k}{\mathcal{Z}^a_{\varepsilon,N}} \left\{ \sum_{\substack{A \subset \{1,\ldots,N-1\} \\ |A|=k}} \int e^{-\mathcal{H}_{[-1,N+1]}(\varphi)} \prod_{m \in A} \delta_0(d\varphi_m) \prod_{n \in A^{\complement}} d\varphi_n \mathbf{1}_{(\varphi_n \in \Omega^a)} \right\}.
\end{aligned}$$

The term in braces in the r.h.s. is a positive number depending on $a, k, N$ but not on $\varepsilon$: let us call it $C^a(k,N)$. Summing over $k = 0, \ldots, N-1$ we obtain

$$\mathcal{Z}^a_{\varepsilon,N} = \sum_{k=0}^{N-1} \varepsilon^k C^a(k,N), \qquad \widetilde{\mathsf{F}}^a_N(t) = \frac{1}{N} \log\left( \sum_{k=0}^{N-1} e^{tk} C^a(k,N) \right),$$

where $\widetilde{\mathsf{F}}^a_N(t) := \mathsf{F}^a_N(e^t)$; cf. (1.10). Differentiating it the variable $t$ we have

$$(\widetilde{\mathsf{F}}^a_N)'(t) = \frac{1}{N} \mathbb{E}^a_{e^t,N}(\ell_N), \qquad (\widetilde{\mathsf{F}}^a_N)''(t) = \frac{1}{N} \mathrm{var}_{\mathbb{P}^a_{e^t,N}}(\ell_N) \geq 0,$$

which proves (1.13) and the convexity of $\widetilde{\mathsf{F}}^a_N(t)$.

Now fix $x \in [0,1]$. For every $\alpha \geq 0$ we have

$$\mathbb{P}^a_{\varepsilon,N}(\ell_N/N > x) = \mathbb{P}^a_{\varepsilon,N}(e^{\alpha \ell_N} > e^{\alpha x N}) \leq e^{-\alpha x N}\, \mathbb{E}^a_{\varepsilon,N}(e^{\alpha \ell_N}). \tag{A.1}$$

Using the above relations we can write

$$\begin{aligned}
\mathbb{E}^a_{\varepsilon,N}(e^{\alpha \ell_N}) &= \sum_{k=0}^{N-1} e^{\alpha(k+1)} \mathbb{P}^a_{\varepsilon,N}(\ell_N = k+1) \\
&= e^{\alpha} \sum_{k=0}^{N-1} (e^{\alpha}\varepsilon)^k \frac{C^a(k,N)}{\mathcal{Z}^a_{\varepsilon,N}} = e^{\alpha} \frac{\mathcal{Z}^a_{e^{\alpha}\varepsilon,N}}{\mathcal{Z}^a_{\varepsilon,N}},
\end{aligned}$$

therefore by (1.8) we have $N^{-1} \log \mathbb{E}^a_{\varepsilon,N}(e^{\alpha \ell_N}) \longrightarrow \mathsf{F}^a(e^{\alpha}\varepsilon) - \mathrm{F}^a(\varepsilon)$ as $N \to \infty$. If $\varepsilon \neq \varepsilon^a_c$, the free energy $\mathsf{F}^a$ is differentiable at $\varepsilon$ by Theorem 1.2, therefore as $\alpha \downarrow 0$ we have $\mathrm{F}^a(e^{\alpha}\varepsilon) - \mathrm{F}^a(\varepsilon) = \mathrm{D}^a(\varepsilon) \cdot \alpha + o(\alpha)$, where $\mathrm{D}^a(\varepsilon) = \varepsilon \cdot (\mathrm{F}^a)'(\varepsilon)$. Plugging $x = \mathrm{D}^a(\varepsilon) + \delta$ (with $\delta > 0$) and $\alpha$ small into (A.1) we obtain

$$\mathbb{P}^a_{\varepsilon,N}(\ell_N/N > \mathrm{D}^a(\varepsilon) + \delta) \leq e^{-(const.)N}.$$

With almost identical arguments one shows that $\mathbb{P}^a_{\varepsilon,N}(\ell_N/N < \mathrm{D}^a(\varepsilon) - \delta) \leq e^{-(const.')N}$, therefore (1.14) and (1.15) are proven.

## APPENDIX B: LLT FOR THE INTEGRATED RANDOM WALK

We are going to prove Proposition 2.3. We will stick for conciseness to the first relation in (2.10), the second one being analogous. We recall that the density of the random vector $(B_1, I_1)$ is given by (2.7). Then its characteristic function $\Psi(s,t)$ is given by

$$\Psi(s,t) = \exp(-s^2/2 - t^2/6 - st/2).$$

We denote by $\psi_n(u,v) := \mathbf{E}^{(0,0)}[\exp(i(uY_n + vZ_n))]$ the characteristic function of $(Y_n, Z_n)$. An application of the Fourier-transform inversion formula gives

$$\begin{aligned}&|\sigma^2 n^2 \varphi_n^{(0,0)}(\sigma\sqrt{n}y, \sigma n^{3/2} z) - g(y,z)| \\ &\quad \le \frac{1}{2\pi} \int_{-\infty}^{+\infty} \left| \psi_n\left(\frac{s}{\sigma\sqrt{n}}, \frac{t}{\sigma n^{3/2}}\right) - \Psi(s,t)\right| ds\, dt.\end{aligned}$$

The proof consists in showing that the r.h.s. vanishes as $n \to \infty$. More precisely, following the proof of Theorem 2 in [12], Section XV.5, we consider separately the three domains

$$\begin{aligned}\mathcal{D}_1 &= \{(s^2 + t^2) \le A\},\\ \mathcal{D}_2 &= \{A < (s^2 + t^2) \le B^2 n\},\\ \mathcal{D}_3 &= \{(s^2 + t^2) > B^2 n\},\end{aligned}$$

and we show that, for a suitable choice of the positive constants $A$ and $B$ and for large $n$, the integral in the r.h.s. above is less than $\varepsilon$ on each domain, for every fixed $\varepsilon > 0$.

*The domain $\mathcal{D}_1$.* Denoting by $\xi(u) := \mathbf{E}[\exp(iuX_1)]$ the characteristic function of $X_1$, from (2.1) and (2.2) we have

$$\psi_n(u,v) = \mathbf{E}^{(0,0)}\left[\prod_{m=1}^n e^{i(u+mv)X_{n+1-m}}\right] = \prod_{m=1}^n \xi(u+mv). \tag{B.1}$$

Since by hypothesis $\mathbf{E}(X_1) = 0$ and $\mathbf{E}({X_1}^2) = \sigma^2 \in (0,\infty)$, it follows that

$$\xi(u) = 1 - \frac{\sigma^2}{2}u^2 + o(u^2) \qquad (u \to 0), \tag{B.2}$$

hence, *uniformly for* $(s,t)$ *such that* $(s^2+t^2) \le A$, we get from (B.1)

$$\psi_n\left(\frac{s}{\sigma\sqrt{n}}, \frac{t}{\sigma n^{3/2}}\right) = \exp\left[-\frac{1}{2}\sum_{m=1}^n \left(\frac{s}{\sqrt{n}} + m\frac{t}{n^{3/2}}\right)^2 + o(1)\right] \overset{n\to\infty}{\longrightarrow} \Psi(s,t).$$

Therefore, for every choice of the parameter $A$, we can find $n_0 = n_0(A)$ such that the integral $\int_{\mathcal{D}_1} |\psi_n(\frac{s}{\sigma\sqrt{n}}, \frac{t}{\sigma n^{3/2}}) - \Psi(s,t)|\, ds\, dt$ is smaller than $\varepsilon$ for $n \ge n_0$.

*The domain* $\mathcal{D}_2$. From B.2 it follows that $|\xi(u)| = 1 - \sigma^2 u^2/2 + o(u^2)$ and therefore we can fix $B > 0$ such that

$$|\xi(u)| \le \exp\Big(-\frac{\sigma^2}{4}u^2\Big) \qquad \text{for } |u| \le \frac{2B}{\sigma}.$$

Using (B.1) and some rough bounds, we get for $(s,t) \in \mathcal{D}_2$ and for all $n \in \mathbb{N}$

$$\begin{aligned}\Big|\psi_n\Big(\frac{s}{\sigma\sqrt{n}}, \frac{t}{\sigma n^{3/2}}\Big)\Big| &\le \exp\Bigg[-\frac{1}{4}\sum_{m=1}^{n}\Big(\frac{s}{\sqrt{n}} + m\frac{t}{n^{3/2}}\Big)^2\Bigg] \\ &\le \exp\Bigg[-\frac{1}{4}\Big(s^2 + \frac{t^2}{3} + st\Big)\Bigg].\end{aligned}$$

Then by the triangle inequality

$$\begin{aligned}&\int_{\mathcal{D}_2}\Big|\psi_n\Big(\frac{s}{\sigma\sqrt{n}}, \frac{t}{\sigma n^{3/2}}\Big) - \Psi(s,t)\Big|\, ds\, dt \\ &\qquad \le \int_{\{s^2+t^2>A\}} (e^{-1/4(s^2+t^2/3+st)} + \Psi(s,t))\, ds\, dt,\end{aligned}$$

and note that the r.h.s. can be made smaller than $\varepsilon$ by choosing $A$ large (this fixes $A$).

*The domain* $\mathcal{D}_3$. By the triangle inequality we have

$$\begin{aligned}&\int_{\mathcal{D}_3}\Big|\psi_n\Big(\frac{s}{\sigma\sqrt{n}}, \frac{t}{\sigma n^{3/2}}\Big) - \Psi(s,t)\Big|\, ds\, dt \\ &\quad \le \int_{\{s^2+t^2>B^2 n\}}\Big|\psi_n\Big(\frac{s}{\sigma\sqrt{n}}, \frac{t}{\sigma n^{3/2}}\Big)\Big|\, ds\, dt + \int_{\{s^2+t^2>B^2 n\}} \Psi(s,t)\, ds\, dt.\end{aligned}$$

It is clear that the second integral in the r.h.s. vanishes as $n \to \infty$ and it remains to show that the same is true for the first integral $I_1$. With the change of variables $s/(\sigma\sqrt{n}) = r\cos(\theta)$, $t/(\sigma\sqrt{n}) = r\sin(\theta)$ and using (B.1), we can rewrite $I_1$ as

$$\begin{aligned}I_1 &= \sigma^2 n \int_{\{\theta\in[0,2\pi),\, r>B\}}\Big|\psi_n\Big(r\cos(\theta), \frac{r\sin(\theta)}{n}\Big)\Big|\, r\, dr\, d\theta \\ &= \sigma^2 n \int_{\{\theta\in[0,2\pi),r>B\}} \Bigg\{\prod_{m=1}^{n} |\xi|\Big(r\,\Big(\cos(\theta) + \frac{m}{n}\sin(\theta)\Big)\Big)\Bigg\} r\, dr\, d\theta,\end{aligned} \tag{B.3}$$

where by $|\xi|(\cdot)$ we mean the function $u \mapsto |\xi(u)|$. It is convenient to divide the domain of integration over $\theta$ in the two subsets

$$\Theta_a := \{\theta \in [0,2\pi) : |\cos(\theta)| > 1/2\}, \qquad \Theta_b := \{\theta \in [0,2\pi) : |\cos(\theta)| \le 1/2\},$$

and to split accordingly the integral $I_1 = I_{1,a} + I_{1,b}$, with obvious notation. We are going to show that both $I_{1,a}$ and $I_{1,b}$ vanish as $n \to \infty$.

Since $\xi(\cdot)$ is the characteristic function of the absolutely continuous random variables $X_1$, we have $|\xi(u)| < 1$ for all $u \neq 0$; cf. Lemma 4 in [12], Section XV.1, and moreover $|\xi(u)| \to 0$ as $u \to \infty$ by the Riemann–Lebesgue lemma; cf. Lemma 3 in [12], Section XV.4. Therefore

$$\Delta := \sup_{\{u\in\mathbb{R}:|u|\geq B/10\}} |\xi(u)| < 1.$$

We are ready to bound $I_{1,a}$ and $I_{1,b}$. For $r > B$, $\theta \in \Theta_a$ and for $m \leq \lfloor n/4 \rfloor$ we have

$$r\left|\cos(\theta) + \frac{m}{n}\sin(\theta)\right| \geq \frac{B}{4} \geq \frac{B}{10},$$

and therefore $|\xi|(r(\cos(\theta) + \frac{m}{n}\sin(\theta))) \leq \Delta$. Since $|\xi| \leq 1$, coming back to (B.3) we can bound $I_{1,a}$ from above (for $n \geq 4$) by

$$\begin{aligned} I_{1,a} &\leq \sigma^2 n \Delta^{\lfloor n/4 \rfloor} \int_{\{\theta\in\Theta_a, r>B\}} \left\{ \prod_{m=n-3}^{n} |\xi|\left( r\left(\cos(\theta) + \frac{m}{n}\sin(\theta)\right)\right)\right\} r\,dr\,d\theta \\ &\leq \sigma^2 n \Delta^{\lfloor n/4 \rfloor} \int_{\mathbb{R}^2} \prod_{m=n-3}^{n} |\xi|\left(x + \frac{m}{n}y\right) dx\,dy. \end{aligned} \tag{B.4}$$

The bound for $I_{1,b}$ is analogous: for $r > B$, $\theta \in \Theta_b$ and for $m \geq \lfloor (3n)/4 \rfloor$ we have

$$r\left|\cos(\theta) + \frac{m}{n}\sin(\theta)\right| \geq B\left(\frac{3}{4}\frac{\sqrt{3}}{2} - \frac{1}{2}\right) \geq \frac{B}{10},$$

and therefore $|\xi|(r(\cos(\theta) + \frac{m}{n}\sin(\theta))) \leq \Delta$. Since $|\xi| \leq 1$, in analogy to (B.4) we can bound $I_{1,b}$ from above by

$$\leq \sigma^2 n \Delta^{\lfloor n/4 \rfloor - 4} \int_{\mathbb{R}^2} \prod_{m=n-3}^{n} |\xi|\left(x + \frac{m}{n}y\right) dx\,dy. \tag{B.5}$$

Combining (B.4) and (B.5), we can finally bound $I_1 = I_{1,a} + I_{1,b}$ from above by

$$I_1 \leq 2\sigma^2 n \Delta^{\lfloor n/4 \rfloor - 4} \int_{\mathbb{R}^2} \prod_{m=n-3}^{n} |\xi|\left(x + \frac{m}{n}y\right) dx\,dy. \tag{B.6}$$

Since $\Delta < 1$, if we prove that the integral in the r.h.s. is bounded by $C \cdot n$ for some positive constant $C$, then it follows that $I_1 \to 0$ as $n \to \infty$ and the proof is completed.

Notice that $|\xi|^2(\cdot) = \xi(\cdot)\,\xi^*(\cdot)$ is the characteristic function of the random variable $X_1 - X_2$, which has an absolutely continuous law with *bounded density* [this is because the density of $X_1$, that is $\exp(-V(\cdot))$, is bounded by hypothesis]. Since $|\xi|^2(\cdot) \geq 0$, it follows from the corollary to Theorem 3 in [12], Section XV.3, that $|\xi|^2(\cdot)$ is integrable over the whole real line, that is, $\|\xi\|_2^2 := \int_{\mathbb{R}} |\xi|^2(x)\,dx < \infty$. By Young's inequality, we can bound the integral in the r.h.s. of (B.6) by

$$\int_{\mathbb{R}^2} \prod_{m=n-3}^{n} |\xi|\left(x+\frac{m}{n}y\right) dx\,dy$$
$$\le \frac{1}{2}\int_{\mathbb{R}^2} |\xi|^2(x+y)\cdot|\xi|^2\left(x+\frac{n-1}{n}y\right) dx\,dy$$
$$+\frac{1}{2}\int_{\mathbb{R}^2} |\xi|^2\left(x+\frac{n-2}{n}y\right)\cdot|\xi|^2\left(x+\frac{n-3}{n}y\right) dx\,dy.$$

However, by a simple change of variables it is easy to see that both the integrals in the r.h.s. equal $n\cdot(\|\xi\|_2^2)^2$, and the proof is completed.

## APPENDIX C: ENTROPIC REPULSION

We are going to prove Proposition 1.5. Notice that the lower bound in (1.23) is easy: $\mathbf{P}(\Omega_N^+) \ge \mathbf{P}(Y_1 \ge 0,\dots,Y_N \ge 0) \sim (const.)/\sqrt{N}$, where the last asymptotic behavior is a classical result of fluctuation theory for random walks with zero mean and finite variance; cf. [12]. Moreover the upper bound in (1.24) follows immediately from the upper bound in (1.23) and Lemma 2.4 [recall the definition (2.11)]. Therefore it remains to prove the lower bound in (1.24), or equivalently (2.13), and the upper bound in (1.23).

**C.1. Proof of (2.13).** We want to get a polynomial lower bound for $w_{0,0}(N)$ as $N\to\infty$. The difficulty comes from the fact that the process $\{Z_n\}_n$ is conditioned to come back to zero and therefore the comparison with the process $\{Y_n\}_n$ is not straightforward.

For simplicity we limit ourselves to the odd case $N = 2n+1$ with $n \in \mathbb{N}$. Recalling the definition $\Omega_k^+ := \{Z_1 \ge 0,\dots,Z_k \ge 0\}$, by Lemma 2.1 we can write

$$w_{0,0}(2n+1) = \mathbf{P}^{(0,0)}(\Omega_{2n-1}^+ \,|\, Z_{2n}=0, Z_{2n+1}=0)$$
$$= \int_{(\mathbb{R}^+)^{2n-1}} e^{-\sum_{k=0}^{2n} V(\Delta\varphi_k)} \prod_{k=1}^{2n-1} d\varphi_k,$$

where we fix $\varphi_{-1}=\varphi_0=\varphi_{2n}=\varphi_{2n+1}=0$. The first step is to restrict the integration on the set $\mathcal{C}_n(\varepsilon) := (\mathbb{R}^+)^{2n-1} \cap \{|\varphi_n-\varphi_{n-1}|\le\varepsilon, |\varphi_n-\varphi_{n+1}|\le\varepsilon\}$, on which $|\Delta\varphi_n| \le 2\varepsilon$. Since $V(\cdot)$ is continuous and $V(0)<\infty$, we can choose $\varepsilon$ sufficiently small such that $V(x)\le V(0)+1$ for all $|x|\le 2\varepsilon$, so that in particular $V(\Delta\varphi_n)\le V(0)+1$ on the event $\mathcal{C}_n(\varepsilon)$. This observation yields the lower bound

$$w_{0,0}(2n+1) \ge e^{-(V(0)+1)} \int_{\mathcal{C}_n(\varepsilon)} e^{-\sum_{k=0}^{n-1} V(\Delta\varphi_k)} \cdot e^{-\sum_{n+1}^{2n} V(\Delta\varphi_k)} \prod_{k=1}^{2n-1} d\varphi_k.$$

Setting $\mathcal{C}'_n(\varepsilon) := (\mathbb{R}^+)^{n-1} \cap \{|\varphi_n - \varphi_{n-1}| \le \varepsilon\}$, the symmetry $k \to 2n-k$ gives

$$w_{0,0}(2n+1) \ge e^{-(V(0)+1)} \int_0^\infty d\varphi_n \left( \int_{\mathcal{C}'_n(\varepsilon)} e^{-\sum_{k=0}^{n-1} V(\Delta\varphi_k)} \prod_{k=1}^{n-1} d\varphi_k \right)^2.$$

Restricting the first integration on $[0, n^{5/2}]$ and applying Jensen's inequality we get

$$\begin{aligned} w_{0,0}(2n+1) &\ge \frac{e^{-(V(0)+1)}}{n^{5/2}} \left( \int_0^{n^{5/2}} d\varphi_n \int_{\mathcal{C}'_n(\varepsilon)} e^{-\sum_{k=0}^{n-1} V(\Delta\varphi_k)} \prod_{k=1}^{n-1} d\varphi_k \right)^2 \\ &= \frac{e^{-(V(0)+1)}}{n^{5/2}} \{\mathbf{P}(\Omega_n^+, Z_n \le n^{5/2}, |Z_n - Z_{n-1}| \le \varepsilon)\}^2. \end{aligned} \tag{C.1}$$

We are thus left with giving a polynomial lower bound for the probability appearing inside the braces. We observe that by definition $Z_n - Z_{n-1} = Y_n$ and that we have the inclusion $\Omega_n^+ \supset \{Y_1 \ge 0, \dots, Y_n \ge 0\} =: \Lambda_n^+$. Therefore

$$\begin{aligned} \mathbf{P}(\Omega_n^+, Z_n \le n^{5/2}, |Z_n - Z_{n-1}| \le \varepsilon) &\ge \mathbf{P}(\Lambda_n^+, Z_n \le n^{5/2}, Y_n \le \varepsilon) \\ &\ge \mathbf{P}(\Lambda_n^+, Y_n \le \varepsilon) - \mathbf{P}(Z_n > n^{5/2}). \end{aligned} \tag{C.2}$$

For the second term, Chebyshev's inequality and (2.2) yield

$$\mathbf{P}(Z_n > n^{5/2}) \le \frac{\mathbf{E}(Z_n^2)}{n^5} \le \frac{(const.) n^3}{n^5} = \frac{(const.)}{n^2}. \tag{C.3}$$

For the first term, we are going to use some results from fluctuation theory. We denote by $\{(T_k, H_k)\}_{k \ge 0}$ the weak ascending ladder process associated to the random walk $\{Y_k\}_k$, that is, $(T_0, H_0) = (0,0)$ and $T_{k+1} := \inf\{n > T_k : Y_n \ge Y_{T_k}\}$, $H_k := Y_{T_k}$. The celebrated duality lemma [12], Chapter XII, gives

$$\mathbf{P}(\Lambda_n^+, Y_n \le \varepsilon) = \sum_{k=0}^\infty \mathbf{P}(T_k = n, Y_n \le \varepsilon),$$

and applying Alili and Doney's combinatorial identity [3] we get

$$\begin{aligned} \mathbf{P}(\Lambda_n^+, Y_n \le \varepsilon) &= \sum_{k=1}^\infty \frac{k}{n} \mathbf{P}(H_{k-1} \le Y_n < H_k, Y_n \le \varepsilon) \\ &= \frac{1}{n} \sum_{k=1}^\infty \mathbf{P}(H_{k-1} \le Y_n, Y_n \le \varepsilon). \end{aligned}$$

Considering only the $k = 1$ term in the sum gives

$$\mathbf{P}(\Lambda_n^+, Y_n \le \varepsilon) \ge \frac{1}{n} \mathbf{P}(Y_n \in [0, \varepsilon]) \sim \frac{1}{n} \cdot \frac{(const.)}{\sqrt{n}} \qquad (n \to \infty).$$

Putting together these bounds with (C.1) and (C.2), we have shown that $w_{0,0}(2n+1) \ge (const.)/n^{11/2}$, hence the proof of (2.13) is complete.

**C.2. Proof of the upper bound in 1.23.** Throughout this section $\delta$ will denote a small positive parameter. More precisely, let us set

$$c_+(\delta) := \frac{\log 1/((1/2)+2\delta)}{\log(3/2+\delta)}. \tag{C.4}$$

Since $c_+(0) = \log 2/\log\frac{3}{2} > 1$, by continuity $c_+(\delta) > 1$ for small $\delta > 0$. In the sequel, we are free to fix any positive $\delta$ such that $c_+(\delta) > 1$; in fact we are going to prove that (2.12) holds true with $c_+ = c_+(\delta)$. For convenience, we split the proof into four steps.

*First step.* We introduce the integer-valued sequence $g_n$ defined for $n \in \mathbb{N}$ by

$$g_n := \lfloor \exp\{(\tfrac{3}{2}+\delta)^n\}\rfloor. \tag{C.5}$$

We claim that for every $D > 0$ there exists $n_0 \in \mathbb{N}$ such that for every $n \geq n_0$ we have

$$\mathbf{P}\left(Z_{g_{n+1}-g_n} \geq -D\,(g_n)^{7/4} - (g_{n+1}-g_n)\,D\,(g_n)^{3/4}\right) \leq \tfrac{1}{2}+\delta. \tag{C.6}$$

It is easily checked that for large $n$

$$\tfrac{1}{2}(g_n)^{3/2+\delta} \leq g_{n+1} - g_n \leq 2(g_n)^{3/2+\delta},$$

which implies that

$$\frac{-D(g_n)^{7/4} - (g_{n+1}-g_n)D(g_n)^{3/4}}{(g_{n+1}-g_n)^{3/2}} \longrightarrow 0 \qquad (n\to\infty).$$

We have already observed in Section 2.3 that $Z_k/(\sigma k^{3/2})$ converges in distribution toward $\int_0^1 B_s\,ds$, where $\{B_s\}_s$ is a standard Brownian motion, hence

$$\mathbf{P}(Z_{g_{n+1}-g_n} \geq -D\,(g_n)^{7/4} - (g_{n+1}-g_n)D(g_n)^{3/4}) \stackrel{n\to\infty}{\longrightarrow} P\left(\int_0^1 B_s\,ds \geq 0\right) = \tfrac{1}{2},$$

and (C.6) follows.

*Second step.* We claim that we can fix the positive constant $D$ such that for every $n \in \mathbb{N}$

$$\mathbf{P}(\{Z_{g_n} > D(g_n)^{7/4}\} \cup \{Y_{g_n} > D(g_n)^{3/4}\} | Z_{g_k} \geq 0, \forall k = 1,\dots,n) \leq \delta. \tag{C.7}$$

By Chebyshev's inequality and (2.2) we can write for every $m \in \mathbb{N}$

$$\mathbf{P}(Z_m > Dm^{7/4}) \leq \frac{\mathbf{E}(Z_m^2)}{D^2 m^{7/2}} \leq \frac{(const.)\,m^3}{D^2 m^{7/2}} = \frac{(const.)}{D^2\sqrt{m}},$$

and analogously

$$\mathbf{P}(Y_m > Dm^{3/4}) \leq \frac{\mathbf{E}(Y_m^2)}{D^2 m^{3/2}} \leq \frac{(const.')\,m}{D^2 m^{3/2}} = \frac{(const.')}{D^2\sqrt{m}}.$$

The inclusion bound yields

$$\begin{aligned}\mathbf{P}(Z_{g_k} \geq 0, \forall k = 1, \dots, n) &\geq \mathbf{P}(Y_m \geq 0, \forall m = 1, \dots, g_n) \\ &\geq \frac{(const.'')}{\sqrt{g_n}} \qquad (n \to \infty),\end{aligned}$$

where the last bound is a well-known result of fluctuation theory; cf. [12], Section XII.7. Then

$$\begin{aligned}&\mathbf{P}(\{Z_{g_n} > D(g_n)^{7/4}\} \cup \{Y_{g_n} > D(g_n)^{3/4}\} \,|\, Z_{g_k} \geq 0\,, \forall k = 1, \dots, n) \\ &\quad \leq \frac{\mathbf{P}(Z_{g_n} > D(g_n)^{7/4}) + \mathbf{P}(Y_{g_n} > D(g_n)^{3/4})}{\mathbf{P}(Z_{g_k} \geq 0, \forall k = 1, \dots, n)} \leq \frac{(const.) + (const.')}{(const.'')} \frac{1}{D^2},\end{aligned}$$

and (C.7) follows.

*Third step.* We claim that for every $n \geq n_0$ (defined in the preceding step) the following relation holds:

$$\mathbf{P}(Z_{g_k} \geq 0, \forall k = 1, \dots, n+1) \leq (\tfrac{1}{2} + 2\delta) \cdot \mathbf{P}(Z_{g_k} \geq 0, \forall k = 1, \dots, n). \tag{C.8}$$

Conditioning on the $\sigma$-algebra $\sigma(Z_1, \dots, Z_{g_n}) = \sigma(X_1, \dots, X_{g_n})$ we have

$$\begin{aligned}&\mathbf{P}(Z_{g_k} \geq 0, \forall k = 1, \dots, n+1) \\ &\qquad = \mathbf{E}(\mathbf{P}(Z_{g_{n+1}} \geq 0 | Y_{g_n}, Z_{g_n}) \cdot \mathbf{1}_{\{Z_{g_k} \geq 0, \forall k = 1, \dots, n\}}).\end{aligned} \tag{C.9}$$

Let us introduce the event $\mathcal{A} := \{Z_{g_n} \leq D\,(g_n)^{7/4}, Y_{g_n} \leq D\,(g_n)^{3/4}\}$, the constant $D$ being fixed in the preceding step. By (2.3) we can write

$$\mathbf{P}(Z_{g_{n+1}} \geq 0 | Y_{g_n}, Z_{g_n}) = \mathbf{P}(Z_{g_{n+1} - g_n} \geq -z - (g_{n+1} - g_n) y)|_{y = Y_{g_n}, z = Z_{g_n}},$$

hence on the event $\mathcal{A}$ we have

$$\mathbf{P}(Z_{g_{n+1}} \geq 0 | Y_{g_n}, Z_{g_n}) \leq \mathbf{P}(Z_{g_{n+1}} \geq -D(g_n)^{7/4} - (g_{n+1} - g_n) D(g_n)^{3/4}) \leq \tfrac{1}{2} + \delta,$$

having applied (C.6). Coming back to (C.9), we get

$$\begin{aligned}\mathbf{P}(Z_{g_k} \geq 0, \forall k = 1, \dots, n+1) &\leq (\tfrac{1}{2} + \delta) \mathbf{P}(Z_{g_k} \geq 0, \forall k = 1, \dots, n) \\ &\quad + \mathbf{P}(\mathcal{A}^{\complement}, Z_{g_k} \geq 0, \forall k = 1, \dots, n),\end{aligned}$$

and thanks to (C.7) we can bound the second term in the r.h.s. by

$$\mathbf{P}(\mathcal{A}^{\complement}, Z_{g_k} \geq 0, \forall k = 1, \dots, n) \leq \delta \cdot \mathbf{P}(Z_{g_k} \geq 0, \forall k = 1, \dots, n),$$

hence (C.8) follows.

*Fourth step.* We are finally ready to complete the proof of (2.12). Assume first that $N = g_n$ for some $n \in \mathbb{N}$. Using the inclusion bound and iterating (C.8) we get

$$\mathbf{P}(\Omega_N^+) \le \mathbf{P}(Z_{g_k} \ge 0, \forall k = 1, \dots, n) \le C(\tfrac{1}{2} + 2\delta)^n,$$

where $C > 0$ is an absolute constant. By (C.5) we have

$$g_n = N \quad \Longrightarrow \quad \exp\left\{\left(\frac{3}{2} + \delta\right)^n\right\} \ge N \quad \Longrightarrow \quad n \ge \frac{\log \log N}{\log(3/2 + \delta)},$$

hence

$$\mathbf{P}(\Omega_N^+) \le C\left(\frac{1}{2} + 2\delta\right)^{\log\log N/\log(3/2+\delta)} = \frac{C}{(\log N)^{c_+(\delta)}},$$

where we recall that $c_+(\delta)$ has been defined in (C.4). Since $\delta > 0$ has been chosen sufficiently small such that $c_+(\delta) > 1$, (2.12) is proven. In the general case, let

$$n^* := \max\{n \in \mathbb{N} : g_n \le N\}.$$

Since $g_{k+1} \approx (g_k)^{3/2+\delta}$, it is easily checked that $g_{n^*} \ge \sqrt{N}$ for large $N$ (provided $\delta < \frac{1}{2}$, which is no harm). Therefore we can repeat the above arguments just replacing $n$ by $n^*$ and $N$ by $\sqrt{N}$, thus getting

$$\mathbf{P}(\Omega_N^+) \le \frac{C}{(\log \sqrt{N})^{c_+(\delta)}} = \frac{2^{c_+(\delta)} C}{(\log N)^{c_+(\delta)}},$$

and (2.12) is proven in full generality.

**Acknowledgments.** We are very grateful to Yvan Velenik and Ostap Hryniv for suggesting the present problem to us, and to Giambattista Giacomin and Fabio Lucio Toninelli for enlightening discussions on the content of Section 1.6.

Dipartimento di Matematica Pura e Applicata
Università degli Studi di Padova
via Trieste 63
35121 Padova
Italy
E-mail: francesco.caravenna@math.unipd.it

Technische Universität Berlin
Institut für Mathematik
Strasse des 17. Juni 136
10623 Berlin
Germany
E-mail: deuschel@math.tu-berlin.de